\documentclass[11pt]{article}
\usepackage{amssymb,amsmath,amsthm, amsfonts}

\def\sqr#1#2{{\vcenter{\vbox{\hrule height.#2pt
              \hbox{\vrule width.#2pt height#1pt \kern#1pt \vrule width.#2pt}
          \hrule height.#2pt}}}}
\def\signed #1{{\unskip\nobreak\hfil\penalty50
          \hskip2em\hbox{}\nobreak\hfil#1
          \parfillskip=0pt \finalhyphendemerits=0 \par}}
\def\endpf{\signed {$\sqr69$}}
\def\sqr#1#2{{\vcenter{\vbox{\hrule height.#2pt
              \hbox{\vrule width.#2pt height#1pt \kern#1pt \vrule width.#2pt}
              \hrule height.#2pt}}}}
\def\signed #1{{\unskip\nobreak\hfil\penalty50
              \hskip2em\hbox{}\nobreak\hfil#1
              \parfillskip=0pt \finalhyphendemerits=0 \par}}
\def\endpf{\signed {$\sqr69$}}
\def\3n{\negthinspace \negthinspace \negthinspace }
\def\2n{\negthinspace \negthinspace }
\def\1n{\negthinspace }

\def\={\buildrel \triangle \over =}

\def\O{\Omega}

\def\q{\quad}

\def\limsup{\mathop{\overline{\rm lim}}}
\def\liminf{\mathop{\underline{\rm lim}}}

\def\max{\mathop{\rm max}}
\def\min{\mathop{\rm min}}
\def\exp{\mathop{\rm exp}}
\def\sup{\mathop{\rm sup}}

\def\inf{\hbox{\rm inf$\,$}}

\def\|{\Big |}
\def\({\Big (}
\def\){\Big )}
\def\[{\Big[}
\def\]{\Big]}
\def\be{\begin{equation}}
\def\bel{\begin{equation}\label}
\def\ee{\end{equation}}
\def\bt{\begin{theorem}}
\def\bcd{\begin{condition}}
\def\ecd{\end{condition}}
\def\et{\end{theorem}}
\def\bc{\begin{corollary}}
\def\ec{\end{corollary}}
\def\bde{\begin{definition}}
\def\ede{\end{definition}}
\def\bl{\begin{lemma}}
\def\el{\end{lemma}}
\def\bp{\begin{proposition}}
\def\ep{\end{proposition}}
\def\br{\begin{remark}}
\def\er{\end{remark}}
\def\ba{\begin{array}}
\def\ea{\end{array}}
\def\ed{\end{document}}

\def\square#1{\vbox{\hrule\hbox{\vrule height#1%
     \kern#1\vrule}\hrule}}
\def\rectangle#1#2{\vbox{\hrule\hbox{\vrule height#1%
     \kern#2\vrule}\hrule}}

\font\tenbb=msbm10 \font\sevenbb=msbm7 \font\fivebb=msbm5

\newfam\bbfam
\scriptscriptfont\bbfam=\fivebb \textfont\bbfam=\tenbb
\scriptfont\bbfam=\sevenbb

\newtheorem{lemma}{Lemma}[section]
\newtheorem{remark}{Remark}[section]

\newtheorem{theorem}{Theorem}[section]
\newtheorem{corollary}{Corollary}[section]

\newtheorem{definition}{Definition}[section]
\newtheorem{proposition}{Proposition}[section]
\newtheorem{condition}{Condition}[section]

\makeatletter
   
   \@addtoreset{equation}{section}
\makeatother

\begin{document}

\title{ PROBABILISTIC INTERPRETATION FOR SYSTEMS OF ISAACS EQUATIONS WITH TWO REFLECTING BARRIERS \bf\footnote{Partially
supported by the NSF of P.R.China (No. 10701050; 10671112), Shandong
Province (No. Q2007A04), and National Basic Research Program of
China (973 Program) (No. 2007CB814904)}}

\author{ Rainer Buckdahn\\
{\small D\'{e}partement de Math\'{e}matiques, Universit\'{e} de Bretagne
Occidentale,}\\
 {\small 6, avenue Victor-le-Gorgeu, CS 93837, 29238 Brest
Cedex 3, France.}\\
{\small{\it E-mail: Rainer.Buckdahn@univ-brest.fr.}}\\
 Juan Li\\
{\small Department of Mathematics, Shandong University at Weihai, Weihai 264200, P. R. China.;}\\
{\small Institute of Mathematics, School of Mathematical Sciences,
Fudan University, Shanghai 200433.}\\
{\small {\it E-mail: juanli@sdu.edu.cn.}}\\
} \date{}\maketitle \noindent{\bf Abstract}\hskip4mm In this paper
we investigate zero-sum two-player stochastic differential games
whose cost functionals are given by doubly controlled reflected
backward stochastic differential equations (RBSDEs) with two
barriers. For admissible controls which can depend on the whole past
and so include, in particular, information occurring before the
beginning of the game, the games are interpreted as games of the
type ``admissible strategy" against ``admissible control", and the
associated lower and upper value functions are studied. A priori
random, they are shown to be deterministic, and it is proved that
they are the unique viscosity solutions of the associated upper and
the lower Bellman-Isaacs equations with two barriers, respectively.
For the proofs we make full use of the penalization method for
RBSDEs with one barrier and RBSDEs with two barriers. For this end
we also prove new estimates for RBSDEs with two barriers, which are
sharper than those in~\cite{HH}. Furthermore, we show that the
viscosity solution of the Isaacs equation with two reflecting
barriers not only can be approximated by the viscosity solutions of
penalized Isaacs equations with one barrier,  but also directly by
the viscosity solutions of penalized Isaacs equations without
barrier.
\\

\newpage
\section{\large{Introduction}}

\hskip1cm In this paper we study two-player zero-sum stochastic
differential games whose cost functionals are given by doubly
controlled backward stochastic differential equations (RBSDEs) with
two reflecting barriers. Fleming and Souganidis~\cite{FS1} were the
first to study in a rigorous manner two-player zero-sum stochastic
differential games. They namely proved that the lower and the upper
value functions of such games satisfy the dynamic programming
principle, and that they are the unique viscosity solutions of the
associated Bellman-Isaacs equations and coincide under the Isaacs
condition. So a lot of recent works are based on the ideas developed
in their pioneering paper~\cite{FS1}; see, for
instance,~\cite{BCR},~\cite{BL1},~\cite{BL2},~\cite{HT}. The reader
interested in this subject is also referred to the references given
in~\cite{FS1}.

On the other hand, general non-linear backward stochastic
differential equations (BSDEs) were first introduced by Pardoux and
Peng~\cite{PaPe} in 1990. They have been studied since then by a lot
of authors and have found various applications, namely in stochastic
control, finance and the second order PDE theory (for more details,
see, for instance~\cite{ElPeQu} and the references therein). BSDE
methods, originally developed by Peng~\cite{Pe1},~\cite{Pe2} for the
stochastic control theory, have been introduced in the theory of
stochastic differential games by Hamad\`{e}ne, Lepeltier~\cite{HL1}
and Hamad\`{e}ne, Lepeltier and Peng~\cite{HLP1} for the study of
games with a dynamics whose diffusion coefficient is strictly
elliptic and doesn't depend on controls. In Buckdahn and
Li~\cite{BL1} there isn't any such restriction on the diffusion
coefficient and they used a completely new approach to study
stochastic differential games. In their framework, in difference to
that of~\cite{FS1}, the admissible controls can depend on the whole
past, including information occurring before the beginning of the
game, and, with the help of a Girsanov transformation argument, the
a priori random lower and upper value functions were shown to be
deterministic. This new approach in combination with BSDE methods
(in particular the notion of stochastic backward semigroups, see
Peng~\cite{Pe1}) allowed them to prove the dynamic programming
principle (DPP) for the upper and lower value functions of the game
as well as to study the associated Isaacs equations in a very
straight-forward way (i.e., in particular without making use of so
called $r$-strategies and $\pi$-admissible strategies playing a
crucial role in~\cite{FS1}).

Also in this paper we investigate two-player zero-sum stochastic
differential games, but in difference to the setting chosen in the
papers mentioned above, we consider now more general running cost
functionals which are defined by doubly controlled RBSDEs with two
barriers. Cvitanic and Karatzas~\cite{CK} were the first to
introduce RBSDEs with two reflecting barriers and to study the
existence and the uniqueness for this type of equation. An RBSDE can
be understood as a BSDE whose solution is forced to stay between two
prescribed continuous processes $L$ and $U,$ called the lower and
the upper obstacle, respectively. The forces imposing the reflection
of the first component $Y$ of the solution of the BSDE at the lower
and the upper barrier, respectively, are described by two continuous
increasing processes $K^{+}$ and $K^{-}$; they are a part of the
solution quadruplet of the RBSDE. Cvitanic and Karatzas also
established the connection between RBSDEs and Dynkin games. Their
work has generalized an earlier one by El Karoui, Kapoudjian,
Pardoux, Peng and Quenez~\cite{EKPPQ} studying BSDEs reflected at a
single barrier. Since then, motivated by the various applications of
the RBSDEs, namely those to stochastic differential games as well as
in finance, many authors have worked on this subject. We refer, for
instance, to the
papers~\cite{EPaQ},~\cite{H},~\cite{HL1},~\cite{HL2},~\cite{HLM},
~\cite{HLP1},~\cite{HH},~\cite{MC},~\cite{PX},~\cite{WY} but also to
the references in these papers. In \cite{CK} the authors used two
different approaches for their existence and uniqueness proof, one
is based on Picard's method and Dynkin games, the other on the
penalization method. Both approaches led to different types of
conditions on the barriers: while one approach used the Mokobodski
condition, the other one supposed the regularity in the sense that
the barriers are uniformly approximable by It\^{o} processes.
Hamad\`{e}ne and Hassani~\cite{HH} replaced these conditions, which
are rather difficult to verify in concrete applications, by the
condition that the both barriers $U$ and $L$  are completely
separated, i.e., $L_t<U_t,$ a.s., for all time points $t$. In our
paper we adopt this condition on the barriers.

The cost functionals defined by the doubly controlled RBSDEs with
two barriers are described the payoff for Player I and the cost for
Player II; they are, in particular, random variables. As well known
for stochastic differential games, the players cannot restrict to
play only control processes; one player has to fix a strategy while
the other player chooses the best answer
 to this strategy in form of a control process. So the lower value function $W$ is defined as the essential infimum  of the
 essential supremum of all cost functionals, where the essential infimum is taken over all admissible strategies of Player II and
 the essential supremum is taken over all admissible controls of Player I. The upper value function $U$ is defined by changing the
 roles of the both players: as the essential supremum  of the essential infimum of all cost functionals, where the essential supremum
 is taken over all admissible strategies of Player I and the essential infimum is taken over all admissible controls of Player II;
 for the precise definitions see (3.9) and (3.10). The objective of our paper is to investigate these lower and upper value functions $W$\
and $U$.  The main results of the paper state that $W$ and $U$ are
deterministic (Proposition 3.1) continuous viscosity solutions of
associated Bellman-Isaacs equations with two obstacles (Theorem
4.1), and they satisfy the DPP (Theorem 3.1) .

We emphasize that the random fields $W$ and $U$, introduced as
combination of essential infimum and essential supremum over a class
of essentially bounded random variables, are deterministic is far
from being trivial. For the proof of Proposition 3.1 we adapt the
argument developed by Buckdahn and Li~\cite{BL1}. Proposition 3.1
then allows to prove the DPP in a straight forward way with the help
of the method of stochastic backward semigroups introduced by
Peng~\cite{Pe1}, which is here extended to RBSDEs with two barriers.
Another key element in the proof of the DPP is the improvement of
former estimates for BSDE with two barriers, obtained by~\cite{HH}.
In fact, we prove that, in the Markovian framework, under
appropriate assumptions the dependence of the solution on the random
initial value of the driving SDE (on which also the barrier
processes depend) is Lipschitz continuous (Proposition 6.1).

The proof that the lower and upper value functions are viscosity
solutions of the associated Isaacs equations with obstacles (Theorem
4.1) is based on a penalization method. As a byproduct we obtain
that the viscosity solutions of penalized Isaacs equations with a
lower barrier as well as of those with an upper barrier (see (4.7),
(4.8)) converge to the viscosity solution of the Isaacs equation
with two obstacles (4.1) (Theorem 4.2). We also obtain that the
viscosity solution of penalized Isaacs equations without obstacles
(see (4.14)) converges to the viscosity solution of equation (4.1)
(Theorem 4.3). Moreover, similar to Buckdahn and
Li~\cite{BL1},~\cite{BL2} we prove the uniqueness (Theorem 5.1) in a
class of continuous functions with a growth condition which was
introduced by Barles, Buckdahn and Pardoux~\cite{BBE} and is weaker
than the polynomial growth assumption.

Our paper is organized as follows. Sections 2 recalls some elements
of the theory of BSDES and RBSDEs with one barrier and two barriers,
respectively, which will be needed in the sequel. Section 3
introduces the setting of stochastic differential games with two
reflecting barriers and their lower and upper value functions $W$
and $U$. It is proved there that these both functions are
deterministic (Proposition 3.1) and satisfy the DPP (Theorem 3.1).
In Section 4 we show that $W$ and $U$ are viscosity solutions of the
associated Bellman-Isaacs equations with two barriers (Theorem 4.1);
their uniqueness is studied in Section 5. We also study
approximations of $W$ and $U$ by viscosity solutions of Isaacs
equations with one or even without obstacles (Theorem 4.2 and
Theorem 4.3). Finally, after having characterized $W$ and $U$ as the
unique viscosity solutions of the associated Bellman-Isaacs
equations with two barriers we show that $W$\ is less than or equal
to $U$, and that, under the Isaacs condition, $W$ and $U$ coincide
(one says that the game has a value). For the sake of readability of
the paper the presentation of some basic properties of RBSDEs with
two barriers associated with forward SDEs, which are needed for our
studies, is postponed to the appendix (Section 6). However, it is
not only a recall, some new results on RBSDEs with two barriers are
also given there, namely Proposition 6.1, already mentioned above.
On the other hand, the proof of DPP (Theorem 3.1) is also given in
the appendix. \hskip1cm

\section{ {\large Preliminaries}}

  \hskip1cm The purpose of this section is to introduce some basic notions and
results concerning BSDEs and RBSDEs with one and two barriers, which
will be needed in the subsequent sections. In all that follows we
will work on the classical Wiener space $(\Omega, {\cal{F}}, P)$:
For an arbitrarily fixed time horizon $T>0$, $\Omega$ is the set of
all continuous functions from $[0, T]$ to ${\mathbb{R}}^d$, with
initial value $0$ ($\Omega= C_0([0, T];{\mathbb{R}}^d))$ and $
{\cal{F}} $ is the Borel $\sigma$-algebra over $\Omega$, completed
by the Wiener measure on $P$. On this probability space the
coordinate process $B_s(\omega)=\omega_s,\ s\in [0, T],\ \omega\in
\Omega$, is a $d$-dimensional Brownian motion. By
${\mathbb{F}}=\{{\mathcal{F}}_s,\ 0\leq s \leq T\}$\ we denote the
natural filtration generated by the coordinate process $B$\ and
augmented by all P-null sets, i.e.,
$${\mathcal{F}}_s=\sigma\{B_r, r\leq s\}\vee {\mathcal{N}}_P,\ \  s\in [0, T]. $$
Here $ {\cal{N}}_P$ is the set of all P-null sets.

We also shall introduce the following both spaces of processes which
will be used frequently in the sequel:
$$\begin{array}{lll}
&\cdot {\cal{S}}^2(0, T; {\mathbb{R}}):=\{(\psi_t)_{0\leq t\leq T}\
\mbox{real-valued adapted continuous process}:
E[\sup\limits_{0\leq t\leq T}| \psi_{t} |^2]< +\infty\};\\
&\cdot {\cal{H}}^{2}(0,T;{\mathbb{R}}^{n}):=\{(\psi_t)_{0\leq t\leq
T}\ {\mathbb{R}}^{n}\mbox{-valued progressively measurable
process}:\\
&\ \hskip3cm \parallel\psi\parallel^2=E[\int^T_0|\psi_t|^2dt]<+\infty \};\\
&\cdot {\cal{A}}^2_{\mbox{c}}(0, T;{\mathbb{R}}):=\{(\psi_t)_{0\leq
t\leq T}\ \mbox{real-valued adapted continuous non-decreasing
process
with}\ \psi_0=0:\\
&\ \hskip3cm E[|\psi_{T} |^2]< +\infty \}.
\end{array}
$$
\noindent(Recall that $|z|$ denotes the Euclidean norm of $z\in
    {\mathbb{R}}^{n}$).
 Let us now consider a measurable function $g:
\Omega\times[0,T]\times {\mathbb{R}} \times {\mathbb{R}}^{d}
\rightarrow {\mathbb{R}} $ with the property that $(g(t, y,
z))_{t\in [0, T]}$ is ${\mathbb{F}}$-progressively measurable for
all $(y,z)$ in ${\mathbb{R}} \times {\mathbb{R}}^{d}$. We make the
following standard assumptions on $g $ throughout the paper:
 \vskip0.2cm

(A1) There is some real $C\ge 0$  such that, P-a.s., for all $t\in
[0, T],\
y_{1}, y_{2}\in {\mathbb{R}},\ z_{1}, z_{2}\in {\mathbb{R}}^d,\\
\mbox{ }\hskip4cm   |g(t, y_{1}, z_{1}) - g(t, y_{2}, z_{2})|\leq
C(|y_{1}-y_{2}| + |z_{1}-z_{2}|).$
 \vskip0.2cm

(A2) $g(\cdot,0,0)\in {\cal{H}}^{2}(0,T;{\mathbb{R}})$. \vskip0.2cm

The following result on BSDEs is by now well known; for its proof
the reader is referred, for instance to the pioneering work by
Pardoux and Peng~\cite{PaPe}, but also to~\cite{ElPeQu}.
 \bl Let the function $g$ satisfy the assumptions (A1) and (A2). Then, for any random variable $\xi\in L^2(\O, {\cal{F}}_T,$ $P),$ the
BSDE associated with $(g, \xi)$
 \be Y_t = \xi + \int_t^Tg(s,Y_s,Z_s)ds - \int^T_tZ_s\,
dB_s,\q 0\le t\le T, \label{BSDE} \ee
 has a unique adapted solution
$$(Y_t, Z_t)_{t\in [0, T]}\in {\cal{S}}^2(0, T; {\mathbb{R}})\times
{\cal{H}}^{2}(0,T;{\mathbb{R}}^{d}). $$\el

\noindent Besides the existence and uniqueness result we shall also
recall the comparison theorem for BSDEs (see Theorem 2.2 in El
Karoui, Peng, Quenez~\cite{ElPeQu} or Proposition 2.4 in
Peng~\cite{Pe2}).

\bl (Comparison Theorem) Given two coefficients $g_1$ and $g_2$
satisfying (A1) and (A2) and two terminal values $ \xi_1,\ \xi_2 \in
L^{2}(\Omega, {\cal{F}}_{T}, P)$, we denote by $(Y^1,Z^1)$\ and
$(Y^2,Z^2)$\ the solution of the BSDE with the data $(\xi_1,g_1 )$\
and $(\xi_2,g_2 )$, respectively. Then we have:

{\rm (i) }(Monotonicity) If  $ \xi_1 \geq \xi_2$  and $ g_1 \geq
g_2, \ a.s.$, then $Y^1_t\geq Y^2_t$, for all $t\in [0, T]$, a.s.

{\rm (ii)}(Strict Monotonicity) If, in addition to {\rm (i)}, we
also assume that $P\{\xi_1 > \xi_2\}> 0$, then $P\{Y^1_t>
Y^2_t\}>0,$ for all $\ 0 \leq t \leq T,$\ and in particular, $
Y^1_0> Y^2_0.$ \el

 \vskip0.3cm
 \noindent{\bf\it \textbf{2.1 Reflected BSDEs
with one barrier}} \vskip0.3cm

After this short and very basic recall on BSDEs let us consider now
RBSDEs with one barrier. An RBSDE with one barrier is associated
with a terminal condition $\xi \in L^{2}(\Omega,{\cal{F}}_{T}, P)$,
a generator $g$\ and an ``obstacle" process $\{L_t\}_{0\leq t \leq
T}$. We assume that $\{L_t\}_{0\leq t \leq T}\in {\cal{S}}^2(0,
T;{\mathbb{R}})$ and $L_T \leq \xi, \ \mbox{a.s.}$

 A solution of an RBSDE with one barrier is a triplet $(Y, Z, K)$ of ${\mathbb{F}}$-progressively measurable
processes, taking its values in
$\mathbb{R}\times\mathbb{R}^d\times\mathbb{R}_+$ and
   satisfying the following properties

\medskip

 {\rm (i)} $Y \in {\cal{S}}^2(0, T; {\mathbb{R}}), \, Z \in
 {\cal{H}}^{2}(0,T;{\mathbb{R}}^{d})$\ and $K\in {\cal{A}}^2_{c}(0, T; {\mathbb{R}})$;
\be\mbox{\rm (ii)}  \ Y_t = \xi + \int_t^Tg(s,Y_s,Z_s)ds + K_{T} -
K_{t} - \int^T_tZ_sdB_s,\quad t\in [0,T];\qquad\qquad\qquad\
\label{RBSDE1} \ee

{\rm (iii)} $Y_t \geq L_t$,\ a.s., for any $ t\in [0,T];$
\vskip0.5cm

{\rm (iv)}$ \int_0^T(Y_t - L_t)dK_{t}=0.$ \vskip0.5cm The following
two lemmata are borrowed from Theorem 5.2 and Theorem 4.1,
respectively, of the paper~\cite{EKPPQ}.

 \bl Assume that $g$ satisfies (A1) and (A2), $ \xi \in L^{2}(\Omega, {\cal{F}}_{T}, P)$,
 $\{L_t\}_{0\leq t \leq T}\in {\cal{S}}^2(0, T;{\mathbb{R}})$
, and $L_T \leq \xi\ \ a.s.$ Then RBSDE~(\ref{RBSDE1}) has a unique
solution $(Y, Z, K).$\el

\br For shortness, a given triplet $(\xi, g, L)$ is said to satisfy
the Standard Assumptions if the generator $g$ satisfies (A1) and
(A2), the terminal value $ \xi$\ belongs to $ L^{2}(\Omega,
{\cal{F}}_{T}, P)$, and the obstacle process $L\in {\cal{S}}^2(0,
T;{\mathbb{R}})$\ is such that $L_T \leq \xi, \ \mbox{a.s.}$ \er

\bl (Comparison Theorem) We suppose that two triplets $(\xi_1, g_1,
L^1)$ and $(\xi_2, g_2, L^2)$\ satisfy the Standard Assumptions but
assume only for one of the both coefficients $g_1$\ and $ g_2$\ to
fulfill the Lipschitz condition (A.1). Furthermore, we make the following assumptions:
$$
  \begin{array}{ll}
{\rm(i)}&\xi_1 \leq \xi_2,\ \ a.s.;\\
{\rm(ii)}&g_1(t,y,z) \leq g_2(t,y,z),\ a.s., \hbox{ \it for } (t,y,z)\in [0,T]\times
{\mathbb{R}}\times {\mathbb{R}}^{d};\\
{\rm(iii)}& L_t^1 \leq L^2_t,\ \ a.s., \hbox{ \it for } t\in [0,T]. \\
 \end{array}
  $$
  Let $(Y^1,Z^1, K^1)$ and $(Y^2, Z^2, K^2)$ be adapted solutions of RBSDEs~(\ref{RBSDE1}) with data $(\xi_1, g_1,
  L^1)$ and $(\xi_2, g_2, L^2),$ respectively.  Then, $Y^1_{t} \leq Y^2_{t}, $ for all $t\in [0,T],$ \ a.s.\el

 We will also need the following standard result on RBSDEs with one barrier.

\bl  Let $(Y,Z,K)$ be the solution of the above RBSDE~(\ref{RBSDE1})
with data $(\xi, g, L)$\ satisfying the Standard Assumptions. Then,
there exists a constant $C$\ such that
$$
 E [\sup_{t\leq s\leq T }
|Y_s|^2+\int_t^T|Z_s|^2ds+|K_T-K_t|^2 |{{\cal {F}}_t} ]\leq CE
[\xi^2+\left(\int_t^T g(s,0,0)ds\right)^2+\sup_{t\leq s\leq T} L_s^2
|{{\cal {F}}_t} ]. $$ The constant $C$\ depends only on the
Lipschitz constant of $g$.\el

 Lemma 2.5 is based on Propositions 3.5 in~\cite{EKPPQ}
and its generalization by Proposition 2.1 in~\cite{WY}.

 \vskip0.3cm
 \noindent{\bf\it \textbf{2.2 Reflected BSDEs
with two barriers}} \vskip0.3cm

Let us consider now RBSDEs with two barriers. An RBSDE with two
barriers is associated with a terminal condition $\xi \in
L^{2}(\Omega,{\cal{F}}_{T}, P)$, a generator $g$\ and two barriers
$L:=\{L_t\}_{0\leq t \leq
   T}$\ and $U:=\{U_t\}_{0\leq t \leq
   T}$\ which belong to ${\cal{S}}^2(0, T;{\mathbb{R}})$\ and satisfy $L_t<U_t,\,  0\leq t \leq
   T$, a.s., and $L_T\leq \xi \leq U_T, \ \mbox{a.s.}$\

A solution of an RBSDE with two barriers is a quadruplet $(Y,
Z,K^+,K^{-})$ of ${\mathbb{F}}$-progressively measurable processes,
taking its values in
$\mathbb{R}\times\mathbb{R}^d\times\mathbb{R}_+\times\mathbb{R}_+$
and satisfying the following properties

\medskip

 {\rm (i)} $Y \in {\cal{S}}^2(0, T; {\mathbb{R}}), \, Z \in
 {\cal{H}}^{2}(0,T;{\mathbb{R}}^{d})$\ and $K^+, K^{-} \in {\cal{A}}^2_{c}(0, T; {\mathbb{R}})$;
\be \mbox{\rm (ii)} \ Y_t = \xi + \int_t^Tg(s,Y_s,Z_s)ds + (K^+_{T}
- K^+_{t})-(K^{-}_{T} - K^{-}_{t})-\int^T_tZ_sdB_s,\ t\in
[0,T];\label{RBSDE2}\ee

{\rm (iii)} $L_t \leq Y_t\leq U_t$,\ a.s., for any $ t\in [0,T];$
\vskip0.5cm

{\rm (iv)}
$\int_0^T(Y_t-L_t)dK^+_{t}=\int_0^T(U_t-Y_t)dK^{-}_{t}=0.$

\vskip0.5cm The following two lemmata are borrowed from Theorem 3.7
and Theorem 1.3 in Hamad\`{e}ne and Hassani~\cite{HH}, respectively.

 \bl Assume that $g$ satisfies (A1) and (A2), $ \xi$\ belongs to $L^{2}(\Omega, {\cal{F}}_{T}, P)$,
 and $L,\ U\in {\cal{S}}^2(0, T;{\mathbb{R}})$\ are such that $L_t<U_t,\, 0\leq t \leq
   T$, a.s., and $L_T\leq \xi \leq U_T, \ \mbox{a.s.}$\ Then RBSDE~(\ref{RBSDE2})
has a unique solution $(Y, Z, K^+, K^-).$\el

\br For shortness, a given quadruplet $(\xi, g, L, U)$ is said to
satisfy the Standard Assumptions if the generator $g$ fulfills (A1)
and (A2), the terminal value $ \xi$\ belongs to $ L^{2}(\Omega,
{\cal{F}}_{T}, P)$, and the two barriers $L,\ U$ belong to ${\cal{S}}^2(0,
T;{\mathbb{R}})$\ and satisfy $L_t<U_t,\, 0\leq t \leq
   T$, a.s., and $L_T\leq \xi \leq U_T, \ \mbox{a.s.}$ \er

\bl (Comparison Theorem) We suppose that the both quadruplets
$(\xi_1, g_1, L^1, U^1)$ and $(\xi_2, g_2,$ $ L^2, U^2)$\ satisfy
the Standard Assumptions. Let $(Y^1,Z^1, K^{1+}, K^{1-} )$ and
$(Y^2, Z^2, K^{2+}, K^{2-})$ be adapted solutions of
RBSDE~(\ref{RBSDE2}) with data $(\xi_1, g_1, L^1, U^1)$ and $(\xi_2,
g_2, L^2, U^2),$ respectively. Furthermore, we make the following
assumptions:
$$
  \begin{array}{ll}
{\rm(i)}&\xi_1 \leq \xi_2,\ \ a.s.;\\
{\rm(ii)}&g_1(t,Y_t^2, Z_t^2) \leq g_2(t,Y_t^2, Z_t^2),\ a.s., \mbox{for}\ t\in [0,T];\\
{\rm(iii)}& L_t^1 \leq L^2_t,\ U_t^1 \leq U^2_t, \mbox{for all}\ t\in [0,T], \mbox{ a.s.}\\
 \end{array}
  $$
  Then, $Y^1_{t} \leq Y^2_{t},\  $ for all $t\in [0,T],$ a.s. In addition if:
$$
  \begin{array}{ll}
{\rm(i)}&g_1(t,y, z) \leq g_2(t,y, z),\ a.s., \mbox{for}\ (t,y,z)\in [0,T]\times {\mathbb{R}}\times {\mathbb{R}}^{d};\\
{\rm(ii)}& L_t^1 = L^2_t,\ U_t^1=U^2_t,\ a.s., \mbox{for}\ t\in [0,T], \\
 \end{array}
  $$
then we also have: $K_t^{1-}\leq K^{2-}_t,\ \mbox{and}\ K_t^{1+}\geq
K^{2+}_t,\ \mbox{for all}\ t\in [0,T],$ a.s.\el

  We will also need the following two estimates for RBSDEs with two barriers.

\bl  Let $(Y, Z, K^+, K^{-})$ be the solution of the above
RBSDE~(\ref{RBSDE2}) with data $(\xi, g, L, U)$\ satisfying the
Standard Assumptions. Then, there exists a constant $C$\ such that
$$E [\sup_{t\leq s\leq T}
|Y_s|^2|{{\cal {F}}_t} ] \leq CE [\xi^2+\left(\int_t^T
g(s,0,0)ds\right)^2+\sup_{t\leq s\leq T} L_s^2+\sup_{t\leq s\leq T}
U_s^2 |{{\cal {F}}_t} ],\ t\in [0, T].
$$ The constant $C$\ depends only on the Lipschitz
constant of $g$.\el

\noindent{\it Proof.} Let $(Y^1,Z^1,K^{1,+})$\ be the unique
solution of the RBSDE with one lower reflecting barrier associated
with $(\xi, g, L)$. We notice that this RBSDE with one lower
reflecting barrier constitutes a particular case of
RBSDE~(\ref{RBSDE2}) with data $(\xi, g, L, U^1)$\ if we assume that
$U^1=Y^1\vee U$; indeed, in this case $(Y^1,Z^1,K^{1,+},K^{1,-}=0)$
is the unique solution of RBSDE~(\ref{RBSDE2}) with data $(\xi, g,
L, U^1)$. Then from Lemma 2.7, we have $Y_t\leq Y_t^1,\ t\in [0,
T],$\ a.s. Moreover, from Lemma 2.5, \be
 E[\sup_{t\leq s\leq T} |Y^1_s|^2|{{\cal {F}}_t} ] \leq CE
[\xi^2+\left(\int_t^T g(s,0,0)ds\right)^2+\sup_{t\leq s\leq T}
L_s^2|{{\cal {F}}_t} ],\ t\in [0, T]. \ee Similarly, let
$(Y^2,Z^2,K^{2,-})$\ be the unique solution of the RBSDE with one
upper reflecting barrier, associated with the data $(\xi, g, U)$. We
observe that this RBSDE is a particular case of RBSDE~(\ref{RBSDE2})
with data $(\xi, g, L^2, U)$\ if we assume that $L^2=Y^2\wedge L$;
the unique solution of this RBSDE with two reflecting barriers is
given by $(Y^2,Z^2,K^{2,+}=0,K^{2,-})$. From Lemma 2.7, we have then
$Y^2_t\leq Y_t,\ t\in [0, T],$\ a.s., and from Lemma 2.5 we have \be
 E[\sup_{t\leq s\leq T} |Y^2_s|^2|{{\cal {F}}_t} ] \leq CE
[\xi^2+\left(\int_t^T g(s,0,0)ds\right)^2+\sup_{t\leq s\leq T}
U_s^2|{{\cal {F}}_t} ],\ t\in [0, T]. \ee Finally, from (2.4), (2.5) and
$Y^2_t\leq Y_t\leq Y_t^1,\ t\in [0, T],$ a.s., we get the wished estimate for $Y$.  The proof is
complete.\endpf

 \bl Let
$(\xi,g,L, U)$ and $(\xi^{\prime},g^{\prime}, L, U)$ be two quadruplets
satisfying the above Standard Assumptions. We suppose that
$(Y,Z,K^+, K^-)$\ and $(Y^{\prime},Z^{\prime},K^{'+}, K^{'-})$\ are
the solutions of RBSDE~(\ref{RBSDE2}) with the data $(\xi,g,L,U)$\
and $(\xi^{\prime},g^{\prime},L,U)$, respectively. Then there exists
a constant $C$ such that, with the notations,
$$\Delta\xi=\xi-\xi^{\prime},\qquad \Delta g=g-g^{\prime},\qquad \Delta Y=Y-Y^{\prime};$$
$$\Delta Z=Z-Z^{\prime},\qquad \Delta K^+=K^+-K^{\prime+},\qquad \Delta K^-=K^--K^{\prime-},$$
it holds $$\aligned & E[ \sup_{t\leq s\leq T}|\Delta Y_s|^2
+\int_t^T|\Delta Z_s|^2ds +|\Delta K^+_T-\Delta K^-_T-(\Delta K^+_t-\Delta K^-_t)|^2|{{\cal {F}}_t} ]\\
&\leq C E[|\Delta\xi|^2+ \int_t^T|\Delta g(s,Y_s,Z_s)|^2ds|{{\cal
{F}}_t} ].
\endaligned$$
The constant $C$\ depends only on the Lipschitz constant of $g'$\
and $g$. \el

For the proof the reader is referred to Theorem 2.4 in Peng and Xu~\cite{PX}.

 \br For the Markovian situation where the obstacle process is a
deterministic function, we can improve Lemma 2.9 considerably and
show that $Y$\ is Lipschitz continuous with respect to the possibly
random initial condition of the driving SDE (whose solution
intervenes in the RBSDEs as well as in the obstacles), see
Proposition 6.1 in the Section 6.\er

\section{\large{Stochastic
Differential Games with Two Barriers and Associated Dynamic
Programming Principles }}

\hskip1cm We now introduce the framework for our study of stochastic
differential games with reflection for two players. We denote the
control state space of the first player by $U$, and that of the
second one by $V$; the associated sets of admissible controls will
be denoted by ${\mathcal{U}}$\ and ${\mathcal{V}}$, respectively.
The set ${\mathcal{U}}$\ is formed by all $U$-valued
${\mathbb{F}}$-progressively measurable processes and
${\mathcal{V}}$\ is the set of all $V$-valued
${\mathbb{F}}$-progressively measurable processes. The control state
spaces U and V are supposed to be compact metric spaces.

For given admissible controls $u(\cdot)=(u_s)_{s\in[0,T]}\in {\mathcal{U}}$ and
$v(\cdot)=(v_s)_{s\in[0,T]}\in {\mathcal{V}}$, the according orbit which regards $t$
as the initial time and $\zeta \in L^2 (\Omega ,{\mathcal{F}}_t,
P;{\mathbb{R}}^n)$ as the initial state is defined by the solution
of the following SDE:
  \be
  \left \{
  \begin{array}{llll}
  dX^{t,\zeta ;u, v}_s & = & b(s,X^{t,\zeta; u,v}_s, u_s, v_s) ds + \sigma(s,X^{t,\zeta; u,v}_s, u_s, v_s) dB_s,\ s\in
   [t,T], \\
   X^{t,\zeta ;u, v}_t  & = & \zeta,
   \end{array}
   \right.
  \ee
where the mappings
  $$
  \begin{array}{llll}
  &   b:[0,T]\times {\mathbb{R}}^n\times U\times V \rightarrow {\mathbb{R}}^n \
  \mbox{and}\ \   \sigma: [0,T]\times {\mathbb{R}}^n\times U\times V\rightarrow {\mathbb{R}}^{n\times d} \\
     \end{array}
  $$
  satisfy the following conditions:
  $$
  \begin{array}{ll}
 \rm{(i)}& \mbox{For every fixed}\ x\in {\mathbb{R}}^n,\ b(., x, ., .)\ \mbox{and}\ \sigma(., x, ., .)
    \ \mbox{are continuous in}\ (t,u,v);\\
 \rm{(ii)}&\mbox{There exists a }C>0\ \mbox{such that, for all}\ t\in [0,T],\ x, x'\in {\mathbb{R}}^n,\ u \in U,\ v \in V, \\
   &\hskip1cm |b(t,x,u,v)-b(t,x',u ,v)|+ |\sigma(t,x,u,v)-\sigma(t,x',u, v)|\leq C|x-x'|.\\
  \end{array}
  \eqno{\mbox{(H3.1)}}
  $$

From (H3.1) we can get a global linear growth condition for $b$ and
$\sigma$, i.e., the existence of some $C>0$\ such that, for all $0
\leq t \leq T,\ u\in U,\ v \in V,\  x\in {\mathbb{R}}^n $,
  \be
  |b(t,x,u,v)| +|\sigma (t,x,u,v)| \leq C(1+|x|).
  \ee
As recalled (6.2) in Section 6, it follows that, under the above
assumptions, for any $u(\cdot)\in {\mathcal{U}}$ and $v(\cdot)\in
{\mathcal{V}}$, SDE (3.1) has a unique strong solution. Moreover,
for any $p\geq 2$, there exists $C_p\in \mathbb{R}$\ such that, for
any $t \in [0,T]$, $u(\cdot)\in {\mathcal{U}}, v(\cdot)\in
{\mathcal{V}}$\ and $ \zeta, \zeta'\in L^p (\Omega
,{\mathcal{F}}_t,P;{\mathbb{R}}^n),$\ we also have the following
estimates, P-a.s.:
 \be
\begin{array}{rcl}
E[\sup \limits_{s\in [t,T]}|X^{t,\zeta; u, v}_s -X^{t,\zeta';u,
v}_s|^p|{{\mathcal{F}}_t}]
& \leq & C_p|\zeta -\zeta'|^p, \\
E[ \sup \limits_{s\in [t,T]} |X^{t,\zeta
;u,v}_s|^p|{{\mathcal{F}}_t}] & \leq &
                        C_p(1+|\zeta|^p).
\end{array}
\ee The constant $C_p$ depends only on the Lipschitz and the linear
growth constants of $b$\ and $\sigma$ with respect to $x$.

Let now be given four functions
$$
\Phi: {\mathbb{R}}^n \rightarrow {\mathbb{R}},\ h,\ h': [0, T]\times
{\mathbb{R}}^n \rightarrow {\mathbb{R}},\ f:[0,T]\times
{\mathbb{R}}^n \times {\mathbb{R}} \times {\mathbb{R}}^d \times U
\times V \rightarrow {\mathbb{R}}
$$
that satisfy the following conditions:
$$
\begin{array}{ll}
\rm{(i)}& \mbox{For every fixed}\ (x, y, z)\in {\mathbb{R}}^n \times
{\mathbb{R}} \times {\mathbb{R}}^d , f(., x, y, z,.,.)\
\mbox{is continuous in}\ (t,u,v)\ \mbox{and}\\
&\mbox{there exists a constant}\ C>0 \ \mbox{such that, for all}\
t\in [0,T],\ x, x'\in {\mathbb{R}}^n,\ y, y'\in
{\mathbb{R}},\ z, z'\\
&\in {\mathbb{R}}^d,\ u \in U \ \mbox{and}\ v \in V,\\
&\hskip3cm\begin{array}{l}
|f(t,x,y,z,u,v)-f(t,x',y',z',u,v)| \\
\hskip3cm \leq C(|x-x'|+|y-y'| +|z-z'|);
\end{array}\\
\rm{(ii)}&\mbox{There is a constant}\ C>0 \ \mbox{such that, for
all}\ x, x'\in {\mathbb{R}}^n,\\
 &\mbox{  }\hskip3cm |\Phi (x) -\Phi (x')|\leq C|x-x'|;\\
\end{array}$$
$$
\begin{array}{ll}
 \rm{(iii)}&\mbox{For every fixed}\ x\in {\mathbb{R}}^n, , h(., x),\ h'(., x)\
 \mbox{are continuous in}\ t\ \mbox{and}\ \mbox{there is a constant}\\\
& C>0 \ \ \mbox{such that, for all}\ t\in [0,T],\ x, x'\in {\mathbb{R}}^n,\\
 &\mbox{  }\hskip3cm |h(t, x) -h(t, x')|+|h'(t, x) -h'(t, x')|\leq C|x-x'|.\\
 &\mbox{Moreover,}\\
&\hskip 1cm h(t, x)<h'(t, x),\ h(T, x)\leq \Phi(x)\leq h'(T, x),\
\mbox{for all}\ t\in[0, T],\ x\in
    {\mathbb{R}}^n.\\
 \end{array}
 \eqno {\mbox{(H3.2)}}
 $$
From (H3.2) we see that $f$,\ $h,\ h'$\ and $\Phi$\ also satisfy the
global linear growth condition in $x$, i.e., there exists some
$C>0$\ such that, for all $0 \leq t \leq T,\  u\in U,\ v \in V,\
x\in {\mathbb{R}}^n $,
  \be
   |f(t,x,0,0,u,v)|+|\Phi (x)|+|h(t, x)| +|h'(t, x)|\leq C(1+|x|).
   \ee
Let $t\in [0, T]$. For any $u(\cdot) \in {\mathcal{U}}, $\ $v(\cdot)
\in {\mathcal{V}}$\ and $\zeta \in L^2
(\Omega,{\mathcal{F}}_t,P;{\mathbb{R}}^n)$,  the mappings $\xi:=
\Phi(X^{t,\zeta; u, v}_T)$, $L_s=h(s, X^{t,\zeta; u, v}_s )$,
$U_s=h'(s, X^{t,\zeta; u, v}_s )$ and $g(s,y,z):= f(s,X^{t,\zeta; u,
v}_s,y,z,u_s,v_s)$ satisfy the conditions of Lemma 2.6 on the
interval $[t, T]$. Therefore, there exists a unique solution to the
following RBSDE with two barriers:
$$
\begin{array}{lll}
 &{\rm (i)}Y^{t,\zeta; u, v} \in
{\cal{S}}^2(t, T; {\mathbb{R}}),\ Z^{t,\zeta; u, v} \in
{\cal{H}}^{2}(t,T;{\mathbb{R}}^{d}),\  \mbox{and}\ \
  K^{+,t,\zeta; u, v},\ K^{-,t,\zeta; u, v} \in {\cal{A}}^2_{c}(t, T;
  {\mathbb{R}});\\ \end{array}$$
  \be \begin{array}{lll}
&{\rm (ii)} Y^{t,\zeta; u, v}_s = \Phi(X_T^{t,\zeta; u, v}) +
\int_s^Tf(r,X^{t,\zeta; u, v}_r,Y^{t,\zeta; u, v}_r,Z^{t,\zeta; u,
v}_r,u_r, v_r)dr + (K^{+,t,\zeta; u, v}_{T}-K^{+,t,\zeta; u, v}_{s}) \\
&\ \hskip3cm  - (K^{-,t,\zeta; u, v}_{T}-K^{-,t,\zeta; u, v}_{s})-\int^T_sZ^{t,\zeta; u, v}_rdB_r,\ \ \  s\in [t,T];\ \\
&{\rm(iii)} h(s, X_s^{t,\zeta; u, v})\leq Y^{t,\zeta; u, v}_s \leq
h'(s, X_s^{t,\zeta; u, v}),\ \
\mbox{a.s.},\ \mbox{for any}\ s\in [t,T];\\
&{\rm (iv)}\int_t^T(Y^{t,\zeta; u, v}_r - h(r,X_r^{t,\zeta; u,
v}))dK^{+,t,\zeta; u, v}_{r}=\int_t^T(h'(r,X_r^{t,\zeta; u,
v})-Y^{t,\zeta; u, v}_r)dK^{-,t,\zeta; u, v}_{r}=0,\end{array} \ee
where $X^{t,\zeta; u, v}$\ is introduced by equation (3.1).

 Moreover, from Proposition 6.1, we can see that
there exists some constant $C>0$\ such that, for all $0 \leq t \leq
T,\ \zeta, \zeta' \in L^2(\Omega ,
{\mathcal{F}}_t,P;{\mathbb{R}}^n),\ u(\cdot) \in {\mathcal{U}}\
\mbox{and}\ v(\cdot) \in {\mathcal{V}},$\ P-a.s.,
 \be
\begin{array}{ll}
 {\rm(i)} & |Y^{t,\zeta; u, v}_t -Y^{t,\zeta'; u, v}_t| \leq C|\zeta -\zeta'|; \\
 {\rm(ii)} & |Y^{t,\zeta; u, v}_t| \leq C (1+|\zeta|). \\
\end{array}
\ee

 Now, similar to~\cite{BL1} and~\cite{FS1}, we introduce the following subspaces
 of admissible controls and the definition of admissible strategies for the game:

\noindent\bde\ An admissible control process $u=\{u_r, r\in [t,
s]\}$ (resp., $v=\{v_r, r\in [t, s]\}$) for Player I (resp., II) on
$[t, s]\ (t<s\leq T)$\ is an ${\mathbb{F}}$-progressively
measurable process taking values in U (resp., V). The set of all
admissible controls for Player I (resp., II) on $[t, s]$ is denoted
by\ ${\mathcal{U}}_{t, s}$\ (resp., ${\mathcal{V}}_{t, s}).$\ We
identify two processes $u$\ and $\bar{u}$\ in\ ${\mathcal{U}}_{t,
s}$\ and write $u\equiv \bar{u}\ \mbox{on}\ [t, s],$\ if
$P\{u=\bar{u}\ \mbox{a.e. in}\ [t, s]\}=1.$\ Similarly, we interpret
$v\equiv \bar{v}\ \mbox{on}\ [t, s]$\ for two elements $v$\ and $
\bar{v}$\ of ${\mathcal{V}}_{t, s}$. \ede

\bde A nonanticipative strategy for Player I on $[t, s] (t<s\leq T)$
is a mapping $\alpha: {\mathcal{V}}_{t, s}\longrightarrow
{\mathcal{U}}_{t, s}$ such that, for any ${\mathbb{F}}$-stopping
time $S: \Omega\rightarrow [t, s]$\ and any $ v_1, v_2 \in
{\mathcal{V}}_{t, s}$\ with $ v_1\equiv v_2\ \mbox {on}\
\textbf{[\![}t, S\textbf{]\!]},$ it holds $\alpha(v_1)\equiv
\alpha(v_2)\ \mbox {on}\ \textbf{[\![}t, S\textbf{]\!]}$.\
Nonanticipative strategies for Player II on $[t, s]$, $\beta:
{\mathcal{U}}_{t, s}\longrightarrow {\mathcal{V}}_{t, s}$,  are
defined similarly. The set of all nonanticipative strategies
$\alpha: {\mathcal{V}}_{t,s}\longrightarrow {\mathcal{U}}_{t,s}$ for
Player I on $[t, s]$ is denoted by ${\cal{A}}_{t,s}$. The set of all
nonanticipative strategies $\beta:
{\mathcal{U}}_{t,s}\longrightarrow {\mathcal{V}}_{t,s}$ for Player
II on $[t, s]$ is denoted by ${\cal{B}}_{t,s}$.  (\mbox{Recall
that}\ $\textbf{[\![}t, S\textbf{]\!]}=\{(r,\omega)\in [0, T]\times
\Omega, t\leq r\leq S(\omega)\}.)$\ede

For given control processes $u(\cdot)\in {\mathcal{U}}_{t,T}$\ and
$ v(\cdot)\in{\mathcal{V}}_{t,T} $\ we introduce the following associated cost
functional
 \be
J(t, x; u, v):= Y^{t, x; u, v}_t,\ (t, x)\in [0, T]\times
{\mathbb{R}}^n,\ee for which the process $Y^{t, x; u, v}$ is defined by
RBSDE (3.5).

\noindent From Proposition 6.2 (see: Appendix) we get that, for any
$t\in[0, T]$, $\zeta \in L^2 (\Omega ,{\mathcal{F}}_t ,P;
{\mathbb{R}}^n)$,
 \be J(t, \zeta; u, v) = Y^{t,\zeta; u, v}_t,\
 \mbox{P-a.s.}
\ee We emphasize that $J(t, \zeta; u, v)=J(t, x; u, v)|_{x=\zeta}$\
while $Y^{t,\zeta; u, v}$\ is defined by (3.5).
 Being particularly interested in the case of a deterministic $\zeta$, i.e., $\zeta=x\in {\mathbb{R}}^n$,
 we define the lower value function of our stochastic
differential game with reflection\be W(t,x):= \mbox{essinf}_{\beta
\in {\cal{B}}_{t,T}}\mbox{esssup}_{u \in {\mathcal{U}}_{t,T}}J(t,x;
u,\beta(u)) \ee
 and its upper value function
  \be U(t,x):= \mbox{esssup}_{\alpha \in
{\cal{A}}_{t,T}}\mbox{essinf}_{v \in {\mathcal{V}}_{t,T}}J(t,x;
\alpha(v),v). \ee

The names ``lower value function" and ``upper value function" for
$W$\ and $U$, respectively, are justified later by Remark 5.1.

\br (1) Here the essential infimum and the essential supremum should
be understood as one with respect to indexed families of random
variables (see, e.g., Dunford and Schwartz~\cite{DS},
Dellacherie~\cite{D} or the appendix in Karatzas and
Shreve~\cite{KS2} for detailed discussions). The reader is also referred
to Remark 3.1 in~\cite{BL1}.\\
(2) Obviously, under the assumptions (H3.1)-(H3.2) it is a
consequence of (3.6)-{\rm(ii)} that the lower value function
$W(t,x)$\ as well as the upper value function $U(t,x)$ are
well-defined and essentially bounded, ${\mathcal{F}}_{t}$-measurable
random variables. But it turns out that $W(t,x)$\ and $U(t,x)$\ are
even deterministic. For proving this we adapt the new approach by
Buckdahn and Li~\cite{BL1},~\cite{BL2}. In the sequel we will
concentrate on the study of the properties of $W$. The discussion of
the properties of $U$, which are comparable with those of $W$, can
be carried out in a similar manner.\er
 \bp For any $(t, x)\in [0, T]\times {\mathbb{R}}^n$,
we have $W(t,x)=E[W(t,x)]$, P-a.s. Identifying $W(t,x)$ with its
deterministic version $E[W(t,x)]$\ we can consider $W:[0, T]\times
{\mathbb{R}}^n\longrightarrow {\mathbb{R}}$ as a deterministic
function.\ep

 \noindent \textbf{Proof}. We recall that $\Omega=
C_0([0, T];{\mathbb{R}}^d)$\ and denote by $H$\ the Cameron-Martin
space of all absolutely continuous elements $h\in \Omega$\ whose
derivative  $\dot{h}$\ belongs to $L^2([0, T],{\mathbb{R}}^d).$\ For
any $h \in H$, we define the mapping $\tau_h\omega:=\omega+h,\
\omega\in \Omega. $\ Obviously, $\tau_h: \Omega\rightarrow\Omega$\
is a bijection and its law is given by
$P\circ[\tau_h]^{-1}=\exp\{\int^T_0\dot{h}_sdB_s-\frac{1}{2}\int^T_0|\dot{h}_s|^2ds\}P.$\
Let $(t, x)\in [0, T]\times {\mathbb{R}}^n$\ be arbitrarily fixed
and put $H_t=\{h\in H|h(\cdot)=h(\cdot\wedge t)\}.$\ We split now
the proof in the following steps:
 \vskip0.1cm
\noindent $1^{st}$ step: For any $u\in {\mathcal{U}}_{t,T}, \ v\in
{\mathcal{V}}_{t,T},\ h \in H_t,\ J(t, x; u,v)(\tau_h)= J(t, x;
u(\tau_h),v(\tau_h)),\ \mbox{P-a.s.}$ \vskip0.1cm
 Indeed, for $h \in H_t$\ we apply the Girsanov transformation to SDE (3.1) (with
 $\zeta=x$). Notice that since $h\in H_t,$\ we have $dB_s(\tau_h)=dB_s,\ s\in [t, T]$. We compare
 the thus obtained equation with the SDE
 got from (3.1) by substituting the transformed control
 processes $u(\tau_h), v(\tau_h)$\ for $u$\ and $v$. Then, from the uniqueness of the solution of
(3.1) we get $X_s^{t,x; u,v}(\tau_h)=X_s^{t,x;
u(\tau_h),v(\tau_h)},$ $ \mbox{for any}\ s\in [t, T],\
\mbox{P-a.s.}$\ Furthermore, by a similar Girsanov transformation
argument we get from the uniqueness of the solution of RBSDE (3.5),
$$Y_s^{t,x; u,v}(\tau_h)=Y_s^{t,x; u(\tau_h),v(\tau_h)},\ \mbox{for
any}\ s\in [t, T],\ \mbox{P-a.s.,}$$
$$Z_s^{t,x; u,v}(\tau_h)=Z_s^{t,x; u(\tau_h),v(\tau_h)},\  \mbox{dsdP-a.e. on}\ [t, T]\times\Omega,$$
$$K_s^{+,t,x; u,v}(\tau_h)=K_s^{+,t,x; u(\tau_h),v(\tau_h)},\ \mbox{for
any}\ s\in [t, T],\ \mbox{P-a.s.,}$$
$$K_s^{-,t,x; u,v}(\tau_h)=K_s^{-,t,x; u(\tau_h),v(\tau_h)},\ \mbox{for
any}\ s\in [t, T],\ \mbox{P-a.s.}$$ This implies, in particular,
that
$$J(t, x; u,v)(\tau_h)= J(t, x; u(\tau_h),v(\tau_h)),\
\mbox{P-a.s.}$$
 \vskip0.1cm
\noindent $2^{nd}$ step: For $\beta\in {\cal{B}}_{t,T}, \ h \in
H_t,$\ let $\beta^h(u):=\beta(u(\tau_{-h}))(\tau_h),\ u\in
{\mathcal{U}}_{t,T}.$\ Then $\beta^h\in {\cal{B}}_{t,T}.$
\vskip0.1cm Obviously, $\beta^h$\ maps ${\mathcal{U}}_{t,T}$\ into
${\mathcal{V}}_{t,T}$.\ Moreover, this mapping is nonanticipating.
Indeed, let $S: \Omega\rightarrow [t, T]$\ be an
${\mathbb{F}}$-stopping time and $ u_1, u_2 \in {\mathcal{U}}_{t,
T}$\ with $ u_1\equiv u_2\ \mbox {on}\ \textbf{[\![}t,
S\textbf{]\!]}.$\ Then, obviously, $ u_1(\tau_{-h})\equiv
u_2(\tau_{-h})\ \mbox {on}\ \textbf{[\![}t,
S(\tau_{-h})\textbf{]\!]}$ (notice that $ S(\tau_{-h})\ \mbox{is
still a}$ stopping time), and because $\beta\in {\cal{B}}_{t,T}$\ we
have $\beta(u_1(\tau_{-h}))\equiv \beta(u_2(\tau_{-h}))\ $ $ \mbox
{on}\ \textbf{[\![}t, S(\tau_{-h})\textbf{]\!]}$. Therefore,
$$\beta^h(u_1)=\beta(u_1(\tau_{-h}))(\tau_h)\equiv \beta(u_2(\tau_{-h}))(\tau_h)=\beta^h(u_2)\ \mbox
{on}\ \textbf{[\![}t, S\textbf{]\!]}.$$ \vskip0.1cm

\noindent$3^{rd}$ step: For all $h\in H_t$\ and $\beta\in
{\mathcal{B}}_{t, T}$\ we have:
$$\{\mbox{esssup}_{u \in {\mathcal{U}}_{t,T}}J(t,x;
u,\beta(u))\}(\tau_h)=\mbox{esssup}_{u \in
{\mathcal{U}}_{t,T}}\{J(t,x; u,\beta(u))(\tau_h)\},\ \mbox{P-a.s.}
$$

Indeed, with the notation $I(t,x,\beta):=\mbox{esssup}_{u \in
{\mathcal{U}}_{t,T}}J(t,x; u,\beta(u)),\ \beta\in {\mathcal{B}}_{t,
T},$\ we have \ $I(t,x,\beta)\geq J(t,x; u,\beta(u)),$\ and thus
$I(t,x,\beta)(\tau_h)\geq J(t,x; u,\beta(u))(\tau_h), \
\mbox{P-a.s.,\ for}$ $\mbox{ all}\ u\in {\mathcal{U}}_{t,T}.$\ On
the other hand, for any random variable $\zeta$\ satisfying
$\zeta\geq J(t,x; u,\beta(u))(\tau_h)$\ and hence also
$\zeta(\tau_{-h})\geq J(t,x; u,\beta(u)), \ \mbox{P-a.s.,\ for}\
\mbox{ all}\ u\in {\mathcal{U}}_{t,T},$\ we have\
$\zeta(\tau_{-h})\geq I(t,x,\beta), \ $ $ \mbox{P-a.s.,}$\ i.e.,
$\zeta\geq I(t,x,\beta)(\tau_{h}), \ \mbox{P-a.s.}$\ Consequently,
$$I(t,x,\beta)(\tau_{h})=\mbox{esssup}_{u \in
{\mathcal{U}}_{t,T}}\{J(t,x; u,\beta(u))(\tau_h)\},\ \mbox{P-a.s.}$$
 \vskip0.1cm
 \noindent$4^{th}$ step: $W(t,x)$\ is invariant with respect
 to the Girsanov transformation $\tau_h$, i.e.,
  $$W(t,x)(\tau_{h})=W(t,x), \ \mbox{P-a.s., for any}\ h\in H. $$

Let us first assume that $h\in H_t$. Then, similarly to the third
step we can show that for all $h\in H_t$,
$$\{\mbox{essinf}_{\beta \in
{\mathcal{B}}_{t,T}}I(t,x;\beta)\}(\tau_h)=\mbox{essinf}_{\beta \in
{\mathcal{B}}_{t,T}}\{I(t,x; \beta)(\tau_h)\},\ \mbox{P-a.s.}
$$\ Then, using the results of the former three steps we have, for any $h\in H_t,$
 $$
   \begin{array}{rcl}
   W(t,x)(\tau_{h}) & = & \mbox{essinf}_{\beta \in
{\mathcal{B}}_{t,T}}\mbox{esssup}_{u \in
{\mathcal{U}}_{t,T}}\{J(t,x; u,\beta(u))(\tau_h)\}\\
       & = &  \mbox{essinf}_{\beta \in
{\mathcal{B}}_{t,T}}\mbox{esssup}_{u \in
{\mathcal{U}}_{t,T}}J(t,x; u(\tau_h),\beta^h(u(\tau_h)))\\
& = &  \mbox{essinf}_{\beta \in {\mathcal{B}}_{t,T}}\mbox{esssup}_{u
\in
{\mathcal{U}}_{t,T}}J(t,x; u,\beta^h(u))\\
& = &  \mbox{essinf}_{\beta \in {\mathcal{B}}_{t,T}}\mbox{esssup}_{u
\in
{\mathcal{U}}_{t,T}}J(t,x; u,\beta(u))\\
& = &W(t,x),\ \mbox{P-a.s.,}
   \end{array}
$$
where we have used the relations
$\{u(\tau_h)|u(\cdot)\in{\mathcal{U}}_{t,T}\}={\mathcal{U}}_{t,T},\
\{\beta^h|\beta \in {\mathcal{B}}_{t,T} \}={\mathcal{B}}_{t,T}$\ in
order to obtain the both latter equalities. Therefore,\ for any
$h\in H_t,\ W(t,x)$ $(\tau_{h})= W(t,x),\ \mbox{P-a.s.,}$\ and since
$W(t,x)$\ is ${\mathcal{F}}_{t}$-measurable, we have this relation
even for all $ h\in H.$\ Indeed, recall that our underlying
fundamental space is $\Omega= C_0([0, T];{\mathbb{R}}^d)$\ and that,
due to the definition of the filtration, the ${\cal F}_t$-measurable
random variable $W(t, x)(\omega),\ \omega\in \Omega,$\ depends only
on the restriction of $\omega$\ to the time interval $[0, t]$.

 The result of the $4^{th}$ step
combined with the following auxiliary Lemma 3.1 completes our
proof.\endpf

\bl Let $\zeta$\ be a random variable defined over our classical
Wiener space $(\Omega, {\mathcal{F}}_T, P)$, such that
 $\zeta(\tau_{h})=\zeta,\ \mbox{P-a.s., for any}\ h\in H.$\ Then
$\zeta=E\zeta,\ \mbox{P-a.s.}$\el

For the proof the reader is referred, for instance, to Lemma 3.4 in
Buckdahn and Li~\cite{BL1}.

The first property of the lower value function $W(t,x)$\ which we
present below is an immediate consequence of (3.6) and (3.9).

\bl\mbox{  }There exists a constant $C>0$\ such that, for all $ 0
\leq t \leq T,\ x, x'\in {\mathbb{R}}^n$,\be
\begin{array}{llll}
&{\rm(i)} & |W(t,x)-W(t,x')| \leq C|x-x'|;  \\
&{\rm(ii)} & |W(t,x)| \leq C(1+|x|).
\end{array}
\ee \el

 We now discuss (the generalized) DPP for our stochastic
differential game with reflection (3.1), (3.5) and (3.9). For this
end we have to
 define the family of (backward) semigroups associated with RBSDE
 (3.5). This notion of stochastic backward semigroups was first
 introduced by Peng~\cite{Pe1} and applied to study the DPP for
 stochastic control problems. Our approach adapts Peng's ideas to the framework of stochastic differential games with reflection.

 Given the initial data $(t,x)$, a positive number $\delta\leq T-t$, admissible control
 processes $u(\cdot) \in {\mathcal{U}}_{t, t+\delta},\ v(\cdot) \in {\mathcal{V}}_{t, t+\delta}$\ and a real-valued
 random variable $\eta \in L^2 (\Omega,
{\mathcal{F}}_{t+\delta},P;{\mathbb{R}})$\ such that $h(t+\delta,
X_{t+\delta}^{t,x; u, v})\leq\eta\leq h'(t+\delta,
X_{t+\delta}^{t,x; u, v}),\ \mbox{a.s.}$, we put \be G^{t, x; u,
v}_{s,t+\delta} [\eta]:= \tilde{Y}_s^{t,x; u, v},\ \hskip0.5cm
s\in[t, t+\delta], \ee where $(\tilde{Y}_s^{t,x;u, v},
\tilde{Z}_s^{t,x;u, v}, \tilde{K}_s^{+,t,x;u,
v},\tilde{K}_s^{-,t,x;u, v})_{t\leq s \leq t+\delta}$ is the
solution of the following RBSDE with time horizon $t+\delta$: \be
\begin{array}{lll}
 &{\rm (i)}\, \tilde{Y}^{t,x; u, v}
\in {\cal{S}}^2(t, t+\delta; {\mathbb{R}}), \tilde{Z}^{t,x; u, v}
\in {\cal{H}}^{2}(t,t+\delta;{\mathbb{R}}^{d}),\  \mbox{and}\ \
  \tilde{K}^{+, t,x; u, v}, \tilde{K}^{-, t,x; u, v} \in {\cal{A}}^2_{c}(t, t+\delta; {\mathbb{R}});\\
&{\rm (ii)}\,  \tilde{Y}^{t,x; u, v}_s = \eta +
\int_s^{t+\delta}f(r,X^{t,x; u, v}_r,\tilde{Y}^{t,x; u,
v}_r,\tilde{Z}^{t,x;
u,v}_r,u_r, v_r)dr + (\tilde{K}^{+,t,x; u, v}_{t+\delta}-\tilde{K}^{+,t,x; u, v}_{s}) \\
&\ \hskip3cm - (\tilde{K}^{-,t,x; u, v}_{t+\delta}-\tilde{K}^{-,t,x;
u, v}_{s}) -\int^{t+\delta}_s\tilde{Z}^{t,x; u, v}_rdB_r,\ \ \  s\in [t,t+\delta];\ \\
&{\rm(iii)}\, h(s, X_s^{t,x; u, v})\leq \tilde{Y}^{t,x; u, v}_s \leq
h'(s, X_s^{t,x; u, v}),\ \ \mbox{a.s.},\
\mbox{for any}\ s\in [t,t+\delta];\\
&{\rm (iv)}\,  \int_t^{t+\delta}(\tilde{Y}^{t,x; u, v}_r -
h(r,X_r^{t,x; u, v}))d\tilde{K}^{+,t,x; u,
v}_{r}=\int_t^{t+\delta}(h'(r,X_r^{t,x; u, v})-\tilde{Y}^{t,x; u,
v}_r )d\tilde{K}^{-,t,x; u, v}_{r}=0,\end{array} \ee where $X^{t,x;
u, v}$\ is introduced by equation (3.1).

Then, in particular, for the solution $(Y^{t,x;u, v}, Z^{t,x;u, v},
K^{+,t,x;u, v}, K^{-,t,x;u, v})$\ of the RBSDE with two barriers (3.5)
we have \be G^{t,x;u, v}_{t,T} [\Phi (X^{t,x; u, v}_T)] =G^{t,x;u,
v}_{t,t+\delta} [Y^{t,x;u, v}_{t+\delta}]. \ee
Moreover,
$$
\begin{array}{rcl}
 J(t,x;u, v)& = &Y_t^{t,x;u, v}=G^{t,x;u, v}_{t,T} [\Phi (X^{t,x; u, v}_T)]
  =G^{t,x;u,v}_{t,t+\delta} [Y^{t,x;u, v}_{t+\delta}]\\
  &=&G^{t,x;u,v}_{t,t+\delta} [J(t+\delta,X^{t,x;u, v}_{t+\delta};u, v)],
\end{array}
$$
where the latter equality follows from (3.8) and the relation
$$Y^{t,x,u,v}_{t+\delta}=Y^{t+\delta,X^{t,x,u,v}_{t+\delta},u,v}_{t+\delta},$$
a consequence of the uniqueness of the solution of BSDE (3.5) and
that of the associated forward equation (3.1). In particular, we
have
$$W(t,x) =\mbox{essinf}_{\beta \in {\mathcal{B}}_{t,
T}}\mbox{esssup}_{u \in {\mathcal{U}}_{t,
T}}G^{t,x;u,\beta(u)}_{t,T} [\Phi(X^{t,x;u,\beta(u)}_{T})].$$

\br For the better comprehension of the reader let us point out that
if $f$\ is independent of $(y, z)$\ then, for all $s\in [t,
t+\delta],$
$$G^{t,x;u,v}_{s,t+\delta}[\eta]=E[\eta + \int_s^{t+\delta}
f(r,X^{t,x;u, v}_r,u_{r}, v_{r})dr+ (\tilde{K}^{+,t,x; u,
v}_{t+\delta}-\tilde{K}^{+,t,x; u, v}_{s})-(\tilde{K}^{-,t,x; u,
v}_{t+\delta}-\tilde{K}^{-,t,x; u, v}_{s})|{\cal{F}}_s].$$ \er

 \bt\mbox{}Under the
assumptions (H3.1) and (H3.2), the lower value function $W(t,x)$
obeys the following
 DPP : For any $0\leq t<t+\delta \leq T,\ x\in {\mathbb{R}}^n,$
 \be
W(t,x) =\mbox{essinf}_{\beta \in {\mathcal{B}}_{t,
t+\delta}}\mbox{esssup}_{u \in {\mathcal{U}}_{t,
t+\delta}}G^{t,x;u,\beta(u)}_{t,t+\delta} [W(t+\delta,
X^{t,x;u,\beta(u)}_{t+\delta})].
 \ee
  \et
The proof of Theorem 3.1 is given in the appendix.

\br We emphasize that, unlike Buckdahn and Li~\cite{BL1},
\cite{BL2}, here we won't use DPP to prove that $W$ and $U$\ are the
viscosity solutions of the associated Isaacs with two barriers,
respectively. \er

\section{\large Viscosity Solution of Isaacs Equation with Obstacles: Existence Theorem }

 \hskip1cm In this section we consider the following Isaacs
equations with obstacles \be
 \left \{\begin{array}{ll}
 &\!\!\!\!\! {\rm min}\left\{W(t,x)-h(t,x), {\rm max}\[-\frac{\partial }{\partial t} W(t,x) - H^{-}(t, x, W, DW, D^2W),
 W(t,x)-h'(t,x)\]\right\}=0,\\
 &\!\!\!\!\!  W(T,x) =\Phi (x), \\
 \end{array}\right.
\ee and
 \be
 \left \{\begin{array}{ll}
 &\!\!\!\!\!{\rm min}\left\{U(t,x)-h(t,x), {\rm max}\[-\frac{\partial }{\partial t} U(t,x) - H^{+}(t, x, U, DU,
 D^2U),W(t,x)-h'(t,x)\]\right\}=0,\\
  &\!\!\!\!\!  U(T,x) =\Phi (x),
 \end{array}\right.
\ee associated with the Hamiltonians $$ H^-(t, x, y, q, X)=
\mbox{sup}_{u \in U}\mbox{inf}_{v \in
V}\{\frac{1}{2}tr(\sigma\sigma^{T}X)+ q.b+ f(t, x, y, q.\sigma,
u, v)\}$$ and
$$ H^+(t, x, y, q, X)= \mbox{inf}_{v \in
V}\mbox{sup}_{u \in U}\{\frac{1}{2}tr(\sigma\sigma^{T}X)+ q.b+ f(t, x, y, q.\sigma,
u, v)\},$$ respectively, where $\sigma$ stands for $\sigma(t, x,
 u, v),$ $b$ for $b(t, x,u, v)$, and $t\in [0, T],\ x\in
{\mathbb{R}}^n,\ y\in {\mathbb{R}},\ q\in {\mathbb{R}}^n$
and $X\in {\mathbb{S}}^n$ (Recall that
${\mathbb{S}}^n$ denotes the set of  symmetric $n\times n$-matrices. Here the functions $b, \sigma, f, h,
h'\ \mbox{and}\ \Phi$\ are supposed to satisfy (H3.1) and (H3.2),
respectively.

 In this section we want to prove that
the lower value function $W(t, x)$ introduced by (3.9) is the
 viscosity solution of equation (4.1), while the upper value
function $U(t, x)$ defined by (3.10) is the viscosity solution of
equation (4.2). The uniqueness of the viscosity solution will be
shown in the next section for the class of continuous functions
satisfying some growth assumption which is weaker than the
polynomial growth condition. We first recall the definition of a
viscosity solution of equation (4.1), that for equation (4.2) is
similar. We borrow the definition from Crandall, Ishii and
Lions~\cite{CIL}.

\bde {\rm(i)} A real-valued upper semicontinuous function
$W:[0,T]\times {\mathbb{R}}^n\rightarrow {\mathbb{R}}$ is called a
viscosity subsolution of equation (4.1) if $W(T,x) \leq \Phi (x),
\mbox{for all}\ x \in
  {\mathbb{R}}^n$, and if for all functions $\varphi \in C^3_{l, b}([0,T]\times
  {\mathbb{R}}^n)$ and $(t,x) \in [0,T) \times {\mathbb{R}}^n$ such that $W-\varphi $\ attains a
  local maximum at $(t, x)$, we have
     $$
    \min\left\{W(t,x)-h(t,x), \max\[-\frac{\partial \varphi}{\partial t} (t,x) -  H^{-}(t, x, W, D\varphi, D^2\varphi),
    W(t,x)-h'(t,x)\]\right\}\leq 0; \eqno{(4.1')}
     $$
{\rm(ii)} A real-valued lower semicontinuous function $W:[0,T]\times
{\mathbb{R}}^n\rightarrow {\mathbb{R}}$ is called a viscosity
supersolution of equation (4.1) if $W(T,x) \geq \Phi (x), \mbox{for
all}\ x \in {\mathbb{R}}^n$, and if for all functions $\varphi \in
C^3_{l, b}([0,T]\times
  {\mathbb{R}}^n)$ and $(t,x) \in [0,T) \times {\mathbb{R}}^n$ such that $W-\varphi $\ attains a
  local minimum at $(t, x)$, it holds
     $$
   \min\left\{W(t,x)-h(t,x), \max\[-\frac{\partial \varphi}{\partial t} (t,x)-H^{-}(t, x, W, D\varphi, D^2\varphi),W(t,x)-h'(t,x)\]\right\}\geq 0;\eqno{(4.1'')}
     $$
 {\rm(iii)} A real-valued continuous function $W\in C([0,T]\times {\mathbb{R}}^n )$ is called a viscosity solution of equation (4.1) if it is both a viscosity sub- and a supersolution of equation
     (4.1).\ede
\br \mbox{  }$C^3_{l, b}([0,T]\times {\mathbb{R}}^n)$ denotes the
set of the real-valued functions that are continuously
differentiable up to the third order and whose derivatives of the
orders 1, 2 and 3 are bounded.\er

We now state the main result of this section.

\bt  Under the assumptions (H3.1) and (H3.2) the lower value function $W$\ defined
by (3.9) is a viscosity solution of the Isaacs equation with two
barriers (4.1), while $U$\ defined by (3.10) solves the Isaacs
equation with two barriers (4.2) in the viscosity solution sense.
\et

We will develop the proof of this theorem only for $W$, that of $U$\
is analogous. The proof is mainly based on an approximation of our
RBSDE (3.5) by a sequence of penalized BSDEs with one barrier. This
penalization method for RBSDEs was first studied in [9], Section 6
(pp.719-pp.723). \vskip0.1cm

For each $(t,x)\in [0,T]\times{\mathbb{R}^n}$, and $m\in\mathbf{N}$,
let $(^mY^{t,x;u,v}_s)_{t\leq s\leq T}$\ (respectively,
$(^m\overline{Y}^{t,x;u,v}_s)_{t\leq s\leq T}$)\vskip0.1cm\noindent
be the first component of the unique solution of the BSDE with one
reflecting lower (resp., upper)  \vskip0.1cm\noindent barrier
associated with
$(f(r,X^{t,x;u,v}_r,y,z,u_r,v_r)-m(h'(r,X^{t,x;u,v}_r)-y)^-,\Phi(X^{t,x;u,v}_T),h(r,X^{t,x;u,v}_r))$\
\vskip0.1cm\noindent(respectively,
$(f(r,X^{t,x;u,v}_r,y,z,u_r,v_r)+m(y-h(r,X^{t,x;u,v}_r))^-,\Phi(X^{t,x;u,v}_T),h'(r,X^{t,x;u,v}_r))$)
(recall that the solutions $^mY^{t,x;u,v}$ and
$^m\overline{Y}^{t,x;u,v}$\ exist due to Lemma 2.3). We define \be
J_m(t,x;u,v) := ^mY^{t,x;u,v}_t,\ \ \ u\in{\cal U}_{t, T},\ \ v\in
{\cal V}_{t, T},\ 0\leq t\leq T,\ x\in{\mathbb{R}^n},\ee and
associate the lower value function \be W_m(t,x) :=
\text{essinf}_{\beta\in{\cal B}_{t,T}}\text{esssup}_{u\in{\cal
U}_{t,T}}J_m(t,x;u,\beta(u)),\qquad 0\leq t\leq T,\
x\in{\mathbb{R}^n}.\ee (respectively, \be \overline{J}_m(t,x;u,v) :=
^m\overline{Y}^{t,x;u,v}_t,\ \ \ u\in{\cal U}_{t, T},\ \ v\in {\cal
V}_{t, T},\ 0\leq t\leq T,\ x\in{\mathbb{R}^n},\ee for which we
consider the lower value function \be \overline{W}_m(t,x) :=
\text{essinf}_{\beta\in{\cal B}_{t,T}}\text{esssup}_{u\in{\cal
U}_{t,T}}\overline{J}_m(t,x;u,\beta(u)),\qquad 0\leq t\leq T,\
x\in{\mathbb{R}^n}.)\ee

It is known from Buckdahn and Li~\cite{BL2} that $W_m(t,x)$ defined
in (4.4) is in $C([0,T]\times {\mathbb{R}^n})$, has linear growth in
$x$, and is a continuous viscosity solution of the
following Isaacs equation with one barrier: \be
 \left \{\begin{array}{ll}
 &\!\!\!\!\!\min\{W_m(t,x)-h(t,x), -\frac{\partial}{\partial t}  W_m(t,x) -
 \mbox{sup}_{u \in U}\mbox{inf}_{v \in V}\{\frac{1}{2}\text{tr}(\sigma\sigma^T(t,x,u,v)D^2W_m(t,x))\\
 &+DW_m(t,x).b(t,x,u,v)+f_m(t,x,W_m(t,x),DW_m(t,x).\sigma(t,x,u,v),u,v)\}\}=0,\\
  &\!\!\!\!\!  W_m(T,x) =\Phi (x), \hskip0.5cm x \in {\mathbb{R}}^n,
 \end{array}\right.
\ee where
$$\begin{array}{ll}
&f_m(t,x,y,z,u,v)=f(t,x,y,z,u,v)-m(h'(t,x)-y)^-,\\
 &\ \hskip5cm(t,x,y,z,u,v)\in [0, T]\times {\mathbb{R}}^n\times
{\mathbb{R}}\times {\mathbb{R}}^d\times U\times V. \end{array}$$
Also the function $\overline{W}_m(t,x)$ defined in (4.6) is in
$C([0,T]\times {\mathbb{R}^n})$, has linear growth in $x$, and is a
continuous viscosity solution of the following Isaacs equation with
one barrier: \be
 \left \{\begin{array}{ll}
 &\!\!\!\!\!\max\{\overline{W}_m(t,x)-h'(t,x), -\frac{\partial}{\partial t}  \overline{W}_m(t,x) -
 \mbox{sup}_{u \in U}\mbox{inf}_{v \in V}\{\frac{1}{2}\text{tr}(\sigma\sigma^T(t,x,u,v)D^2\overline{W}_m(t,x))\\
 &+D\overline{W}_m(t,x).b(t,x,u,v)+\overline{f}_m(t,x,\overline{W}_m(t,x),D\overline{W}_m(t,x).\sigma(t,x,u,v),u,v)\}\}=0,\\
  &\!\!\!\!\!  \overline{W}_m(T,x) =\Phi (x), \hskip0.5cm x \in {\mathbb{R}}^n,
 \end{array}\right.
\ee where
$$\begin{array}{ll}
&\overline{f}_m(t,x,y,z,u,v)=f(t,x,y,z,u,v)+m(y-h(t,x))^-,\\
 &\ \hskip5cm(t,x,y,z,u,v)\in [0, T]\times {\mathbb{R}}^n\times
{\mathbb{R}}\times {\mathbb{R}}^d\times U\times V. \end{array}$$

\noindent We have the uniqueness of the viscosity solutions $W_m,\
\overline{W}_m$\ in the space $\tilde{\Theta}$\ which is defined by

$\tilde{\Theta}=\{ \varphi\in C([0, T]\times {\mathbb{R}}^n):
\exists\ \widetilde{A}>0\ \mbox{such
 that}$ \vskip 0.1cm
 $\mbox{ }\hskip2cm \lim_{|x|\rightarrow \infty}\varphi(t, x)\exp\{-\widetilde{A}[\log((|x|^2+1)^{\frac{1}{2}})]^2\}=0,\
 \mbox{uniformly in}\ t\in [0, T]\}.$

\bl For all $(t, x)\in [0, T]\times {\mathbb{R}}^n$\ and all $m\geq
1$, $${W}_1(t, x)\geq \cdots \geq {W}_m(t, x) \geq {W}_{m+1}(t,
x)\geq \cdots \geq W(t, x).$$$$\overline{W}_1(t, x)\leq \cdots \leq
\overline{W}_m(t, x) \leq \overline{W}_{m+1}(t, x)\leq \cdots \leq
W(t, x).$$ \el

\noindent{\bf Proof}. Let $m\geq 1$, since $f_m(t,x,y,z,u,v)\geq
f_{m+1}(t,x,y,z,u,v)$, for all $(t,x,y,z,u,v)$\ we obtain from the
comparison theorem for BSDEs with one barrier (Lemma 2.4) that
 $$J_m(t,x,u,v)={}^mY^{t,x;u,v}_t\geq {}^{m+1}Y^{t,x;u,v}_t=J_{m+1}(t,x,u,v), \ \mbox{P-a.s., for any}
 \ u\in{\cal U}_{t, T}\ \mbox{and}\ v\in {\cal V}_{t, T}.$$
Consequently, $W_m(t, x)\geq W_{m+1}(t, x),\ \mbox{for all}\ (t,
x)\in [0, T]\times {\mathbb{R}}^n,\ m\geq 1.$

From the comparison principle of Section 3 [pp.247-pp.256] in~\cite{HH} we
get that, for each $0\leq t\leq T$, $x\in{\mathbb{R}}^n,\
u\in{\cal U}_{t, T}\ \mbox{and}\ v\in {\cal V}_{t, T}$, \be
J_m(t,x;u,v)\geq J(t,x;u,v),\ \mbox{P-a.s.}\ee It follows that
$W_m(t, x)\geq W(t, x),\ \mbox{for all}\ (t, x)\in [0, T]\times
{\mathbb{R}}^n,\ m\geq 1.$ The proof for $\overline{W}_m(t,
x)$ can be carried out in a similar way.\endpf

\br The above lemma allows to introduce the upper semicontinuous
function $\widetilde{{W}}$\ as limit over the non-increasing
sequence of continuous functions ${W}_m,\ m\geq 1,$\ and we have
$$W_1(t,x)\geq \widetilde{{W}}(t,x)(=\lim_{m\uparrow \infty}\downarrow {W}_m(t, x))\geq W(t, x),\ \ (t, x)\in [0, T]\times
{\mathbb{R}}^n.$$ \er

\br The above lemma also allows to introduce the lower semicontinuous
function $\widetilde{\overline{W}}$\ as limit over the
non-decreasing sequence of continuous functions $\overline{W}_m,\
m\geq 1.$\ From
$$\overline{W}_1(t,x)\leq \widetilde{\overline{W}}(t,x)(=\lim_{m\uparrow \infty}\uparrow \overline{W}_m(t, x))\leq W(t, x),\ \ (t, x)\in [0, T]\times
{\mathbb{R}}^n,$$ \noindent and from Remark 4.2 and the linear
growth of $\overline{W}_1$\ and ${W_1}$\ we conclude that also
$\widetilde{{W}}$\ and $\widetilde{\overline{W}}$\ have at most
linear growth. \er

Our objective is to prove that $\widetilde{W}$,
$\widetilde{\overline{W}}$\ and $W$\ coincide and equation (4.1)
holds in viscosity sense. For this end we first prove the following
proposition:

\bp  Under the assumptions (H3.1) and (H3.2) the function
$\widetilde{{W}}(t,x)$\ is a viscosity subsolution of Isaacs
equation (4.1). \ep

\noindent{\bf Proof}. Let $(t,x)\in [0,T)\times {\mathbb{R}}^n$\ and
let $\varphi\in C^3_{l, b}([0,T]\times {\mathbb{R}}^n)$\ be such
that $\widetilde{{W}}-\varphi<\widetilde{{W}}(t,x)-\varphi(t,x)$\
everywhere on $([0,T]\times {\mathbb{R}}^n) -\{(t,x)\}.$\ Then,
since $\widetilde{{W}}$\ is upper semicontinous and
${W}_m(t,x)\downarrow \widetilde{{W}}(t,x)$, $0\leq t\leq T $,
$x\in{\mathbb{R}}^n$, there exists some sequence $(t_m,x_m),\ m\geq
1,$\ such that, at least along a subsequence,

\smallskip

i)$(t_m,x_m)\rightarrow (t,x)$, as $m\rightarrow +\infty$;

\smallskip
ii) $W_m-\varphi\leq W_m(t_m,x_m)- \varphi(t_m,x_m)$ in a
neighborhood of $(t_m,x_m)$, for all $m\geq 1$;

\smallskip
iii) $W_m(t_m,x_m)\rightarrow \widetilde{W}(t,x)$, as $m\rightarrow
+\infty$.

\medskip

\medskip

\noindent From the definition of $\widetilde{{W}}(t,x)$ we know
$\widetilde{{W}}(t,x)\geq h(t,x)$. Therefore we only need to distinguish
two cases. In the case, for which $\widetilde{{W}}(t,x)= h(t,x)$, equation $(4.1')$ is trivially satisfied, and the proof is complete. Let us discuss the second case: $\widetilde{{W}}(t,x)> h(t,x)$. For this case we get as an immediate consequence of i) and iii) the following result:

\smallskip
iv) There exists $N$\ such that ${W}_m(t_m,x_m)>h(t_m,x_m)$, for all
$m>N$.

\noindent Thus, because $W_m$ is a viscosity solution and
hence a subsolution of equation (4.7), we have, for all $m\geq 1$,
\be\begin{array}{lll} &\frac{\partial}{\partial t}\varphi(t_m,x_m)
+\sup_{u\in U}\inf_{v\in
V}\big\{\displaystyle\frac{1}{2}tr(\sigma\sigma^*(t_m,x_m,u,v)
D^2\varphi(t_m,x_m))\\
& +b(t_m,x_m,u,v)D \varphi(t_m,x_m)+f(t_m,x_m,W_m(t_m,x_m),
D\varphi(t_m,x_m)\sigma(t_m,x_m,u,v),u,v)\big\}\\
& -m(h'(t_m,x_m)-W_m(t_m,x_m))^{-}\\
& \geq 0. \end{array}\ee

\noindent Therefore,

\smallskip

$\frac{\partial}{\partial t}\varphi(t_m,x_m) +\sup_{u\in
U}\inf_{v\in
V}\big\{\displaystyle\frac{1}{2}$tr$(\sigma\sigma^*(t_m,x_m,u,v)
D^2\varphi(t_m,x_m))$

$+b(t_m,x_m,u,v)D \varphi(t_m,x_m)+f(t_m,x_m,W_m(t_m,x_m),
D\varphi(t_m,x_m)\sigma(t_m,x_m,u,v),u,v)\big\}$

$\geq 0.$

\medskip

\noindent We recall that $(t_m,x_m)\rightarrow (t,x)$ and
$W_m(t_m,x_m)\rightarrow \widetilde{W}(t,x)$, as $m\rightarrow
+\infty$. On the other hand, from the continuity of the functions
$b,\sigma$ and $f$ we have, in particular, their uniform continuity
on compacts (recall that $U,V$ are compacts). Consequently,

\medskip

$\displaystyle\frac{\partial}{\partial t}\varphi(t_m,x_m)
+\displaystyle\frac{1}{2}$tr$(\sigma\sigma^*(t_m,x_m,u,v)
D^2\varphi(t_m,x_m))$

$+b(t_m,x_m,u,v)D\varphi(t_m,x_m)+f(t_m,x_m,W_m(t_m,x_m),
D\varphi(t_m,x_m)\sigma(t_m,x_m,u,v),u,v)$

\smallskip

\noindent converges, uniformly in $(u,v)$, towards

$\displaystyle\frac{\partial}{\partial t}\varphi(t,x)
+\displaystyle\frac{1}{2}$tr$(\sigma\sigma^*(t,x,u,v)
D^2\varphi(t,x))$

$+b(t,x,u,v)D \varphi(t,x)+f(t,x,\widetilde{W}(t,x),
D\varphi(t,x)\sigma(t,x,u,v),u,v).$

\smallskip

\noindent Therefore, \be\begin{array}{lll}&\frac{\partial}{\partial
t}\varphi(t,x) +\sup_{u\in U}\inf_{v\in
V}\big\{\displaystyle\frac{1}{2}$tr$(\sigma\sigma^*(t,x,u,v)
D^2\varphi(t,x))\\
&+b(t,x,u,v)D\varphi(t,x)+f(t,x,\widetilde{W}(t,x),
D\varphi(t,x)\sigma(t,x,u,v),u,v)\big\}\\
&\geq 0.\end{array}\ee \noindent The above calculation shows that if
$\widetilde{W}(t,x)\leq h'(t,x)$\ then we can conclude
$\widetilde{W}$ is a viscosity subsolution of (4.1). For proving
that $\widetilde{W}(t,x)\leq h'(t,x)$\ we return to the above
inequality (4.10), from where

\medskip
$m(h'(t_m,x_m)-W_m(t_m,x_m) )^{-}$

$\leq\frac{\partial}{\partial t}\varphi(t_m,x_m) +\sup_{u\in
U}\inf_{v\in
V}\big\{\displaystyle\frac{1}{2}$tr$(\sigma\sigma^*(t_m,x_m,u,v)
D^2\varphi(t_m,x_m))$

$+b(t_m,x_m,u,v)D\varphi(t_m,x_m)+f(t_m,x_m,W_m(t_m,x_m),
D\varphi(t_m,x_m)\sigma(t_m,x_m,u,v),u,v)\big\}.$

\smallskip
\noindent When $m$\ tends to $+\infty$\ the limit of the right-hand
side of the above inequality, given by the left hand side of (4.11),
is a real number. Therefore, the left-hand side of the above
inequality cannot tend to $+\infty$. But this is only possible if
$(h'(t_m,x_m)-W_m(t_m,x_m) )^{-}\rightarrow 0$, i.e., if
$\widetilde{W}(t,x)\leq h'(t,x)$. The proof is complete.
\endpf
\bp  Under the assumptions (H3.1) and (H3.2) the function
$\widetilde{\overline{W}}(t,x)$\ is a viscosity supersolution of
Isaacs equations (4.1). \ep

The proof is similar to that of Proposition 4.1, so we omit it.

\medskip
 \noindent \textbf{Proof of Theorem 4.1}. From Theorem 5.1 which is proved in Section 5, Propositions 4.1 and
 4.2 we get $\widetilde{W}(t,x)\leq \widetilde{\overline{W}}(t,x)$. Furthermore, from Remark 4.2 and Remark 4.3 we get
 $\widetilde{W}(t,x)=\widetilde{\overline{W}}(t,x)=W(t,x)$. The proof is complete.\endpf
\medskip

As a byproduct to the proof of Theorem 5.1 we have that both the
viscosity solution $W_m$\ of the Isaacs equation with one obstacle (4.7) and the viscosity
solution $\overline{W}_m$\ of the Isaacs equation with one obstacle (4.8) converge
pointwise to the viscosity solution of the Isaacs equation with two
obstacles (4.1)  :

  \bt $W_m(t,x)\downarrow W(t,x)$\ and  $\overline{W}_m(t,x)\uparrow W(t,x)$, as $m\rightarrow+\infty,$\ for any $(t,x)\in [0, T]\times
  {\mathbb{R}}^n.$\et

On the other hand, we can also describe $W$ as limit of solutions of
a sequence of Bellman-Isaacs equations without obstacle. For this
end we let, for each $(t,x)\in [0,T]\times{\mathbb{R}^n}$, and $m,
n\in\mathbf{N}$,  \vskip0.1cm\noindent
 $(^{m,n}Y^{t,x;u,v}_s)_{t\leq s\leq T}$\
 be the first component of the unique solution of the
BSDE associated with  \vskip0.1cm\noindent
$(f(r,X^{t,x;u,v}_r,y,z,u_r,v_r)-m(h'(r,X^{t,x;u,v}_r)-y)^-+n(y-h(r,X^{t,x;u,v}_r))^-,$
$\Phi(X^{t,x;u,v}_T))$\ (recall that \vskip0.1cm\noindent due to
Lemma 2.1 $^{m,n}Y^{t,x;u,v}$\ exists). We define \be
J_{m,n}(t,x;u,v) := ^{m,n}Y^{t,x;u,v}_t,\ \ \ u\in{\cal U}_{t, T},\
\ v\in {\cal V}_{t, T},\ 0\leq t\leq T,\ x\in{\mathbb{R}^n},\ee and
consider the lower value function \be W_{m,n}(t,x) :=
\text{essinf}_{\beta\in{\cal B}_{t,T}}\text{esssup}_{u\in{\cal
U}_{t,T}}J_{m,n}(t,x;u,\beta(u)),\qquad 0\leq t\leq T,\
x\in{\mathbb{R}^n}.\ee It is known from Buckdahn and Li~\cite{BL1}
that $W_{m,n}(t,x)$ defined in (4.13) is in $C([0,T]\times
{\mathbb{R}^n})$, has linear growth in $x$, and is the unique
continuous viscosity solution of the following Isaacs equations: \be
 \left \{\begin{array}{ll}
 &\!\!\!\!\!-\frac{\partial}{\partial t} W_{m,n}(t,x) -
 \mbox{sup}_{u \in U}\mbox{inf}_{v \in V}\{\frac{1}{2}\text{tr}(\sigma\sigma^T(t,x,u,v)D^2W_{m,n}(t,x))\\
 &+DW_{m,n}(t,x).b(t,x,u,v)+f_{m,n}(t,x,W_{m,n}(t,x),DW_{m,n}(t,x).\sigma(t,x,u,v),u,v)\}\}=0,\\
  &\!\!\!\!\!  W_{m,n}(T,x) =\Phi (x), \hskip0.5cm x \in {\mathbb{R}}^n,
 \end{array}\right.
\ee where
$$\begin{array}{ll}
&f_{m,n}(t,x,y,z,u,v)=f(t,x,y,z,u,v)-m(h'(t,x)-y)^-+n(y-h(t,x))^-,\\
 &\ \hskip5cm(t,x,y,z,u,v)\in [0, T]\times {\mathbb{R}}^n\times
{\mathbb{R}}\times {\mathbb{R}}^d\times U\times V. \end{array}$$

\noindent Furthermore, from Theorem 4.2 in~\cite{BL2} we know, for
any $m\in\mathbf{N}$, we have \be\lim_{n\rightarrow\infty}\uparrow
W_{m,n}(t,x)=W_{m}(t,x).\ee Similarly,  for any $m\in\mathbf{N}$, we
have \be\lim_{n\rightarrow\infty}\downarrow
W_{n,m}(t,x)=\overline{W}_{m}(t,x). \ee
 \bt $\lim_{m\rightarrow\infty}{W}_{m,m}(t,x)=W(t,x)$
 \ for any $(t,x)\in [0, T]\times
  {\mathbb{R}}^n.$\et
\noindent {\textbf{Proof}.} From the comparison theorem for BSDEs
(Lemma 2.2) and the definition of $W_{m,n}(t,x)$ we can get that for
any $n\geq m,$ $ W_{n,m}(t,x)\leq W_{n,n}(t,x)\leq W_{m,n}(t,x)$.
Combining this with (4.15) and (4.16) we get
$$\overline{W}_{m}(t,x)\leq
\underline{\lim}_{n\rightarrow\infty}W_{n,n}(t,x)\leq
\overline{\lim}_{n\rightarrow\infty}W_{n,n}(t,x)\leq W_{m}(t,x),$$
and then taking the limit as $m\rightarrow\infty$ from Theorem 4.2
we have the wished result.
\endpf
\section{\large Viscosity Solution of Isaacs' Equation with obstacles: Uniqueness Theorem }
\ \hskip 0.5cm The objective of this section is to study the
uniqueness of the viscosity solution of Isaacs' equation (4.1),
\hskip1cm \be
 \left \{\begin{array}{ll}
 &\!\!\!\!\! {\rm min}\left\{W(t,x)-h(t,x),\max\[ -\frac{\partial }{\partial t} W(t,x) - H^{-}(t, x, W, DW,
 D^2W),W(t,x)-h'(t,x)\]\right\}=0,\\
  &\!\!\!\!\!  W(T,x) =\Phi (x).
 \end{array}\right.
\ee  associated with the Hamiltonian $$ H^-(t, x, y, q, X)=
\mbox{sup}_{u \in U}\mbox{inf}_{v \in
V}\{\frac{1}{2}tr(\sigma\sigma^{T}(t, x,
 u, v)X)+ q.b(t, x, u, v)+ f(t, x, y, q.\sigma,
u, v)\},$$
 $ t\in [0, T],\ x\in {\mathbb{R}}^n,\ y\in
{\mathbb{R}},\ q\in {\mathbb{R}}^n,\ X\in {\mathbb{S}}^n$.
 The functions $b, \sigma, f,h ,h'\ \mbox{and}\ \Phi$\ are still supposed to satisfy (H3.1) and (H3.2), respectively.

 For the proof of the uniqueness of the viscosity solution for equation (5.1) in the
 space of functions\\
 $\mbox{ }\hskip1.5cm \Theta=\{ \varphi: [0, T]\times {\mathbb{R}}^n\rightarrow {\mathbb{R}}| \exists\ \widetilde{A}>0\ \mbox{such
 that}$\\
 $\mbox{ }\hskip2.5cm \lim_{|x|\rightarrow \infty}\varphi(t, x)\exp\{-\widetilde{A}[\log((|x|^2+1)^{\frac{1}{2}})]^2\}=0,\
 \mbox{uniformly in}\ t\in [0, T]\}$  \vskip 0.1cm \noindent we borrow the main idea from Barles, Buckdahn,
Pardoux~\cite{BBE}. This growth condition was introduced
in~\cite{BBE} to prove the uniqueness of the viscosity solution of
an integro-partial differential equation associated with a decoupled
FBSDE with jumps. It was shown in~\cite{BBE} that this kind of
growth condition is optimal for the uniqueness and can, in general,
not be weakened, even not for PDEs. We adapt the ideas developed
in~\cite{BBE} and~\cite{BL1},~\cite{BL2} to Isaacs' equation (5.1)
to prove the uniqueness of the viscosity solution in $\Theta$. Since
the proof of the uniqueness in $\Theta$\ for equation (4.2) is the
same we will restrict ourselves to that of (5.1). Before stating the
main result of this section, let us begin with two auxiliary
lemmata. Denoting by $K$\ a Lipschitz constant of $f(t,x,.,.)$, that
is uniform in $(t, x),$\ we have the following

\bl Let an upper semicontinuous function $u_1 \in \Theta$\ be a
viscosity subsolution and a lower semicontinuous function $u_2 \in
\Theta$\ be a viscosity supersolution of equation (5.1). Then, the
upper semicontinuous function $\omega:= u_1-u_2$\ is a viscosity
subsolution of the equation
    \be\left\{
         \begin{array}{lll}
     &\!\!\!\!\!{\rm min}\{\omega(t,x),-\frac{\partial }{\partial t} \omega(t,x)
     - {\rm sup}_{u \in
U, v \in V}( \frac{1}{2}tr(\sigma\sigma^{T}(t, x,
 u, v)D^2\omega)+ D\omega.b(t, x, u, v)+ K|\omega|\\
 &\!\!\!\!\!\mbox{ }\hskip1cm  +K|D\omega .\sigma(t, x, u, v)|) \}= 0, \ \hskip2cm  (t, x)\in [0, T)\times
 {\mathbb{R}}^n,\\
&\!\!\!\!\!\omega(T,x) =0,\ \hskip1cm  x \in {\mathbb{R}}^n.
     \end{array}
\right.
   \ee
  \el
\noindent{\bf Proof}. The proof is similar to that of Lemma 3.7
in~\cite{BBE} or Lemma 5.1 in~\cite{BL1}, the main difference
consists in the fact that here we have to deal with an obstacle
problem.

We observe that $\omega(T,x)=u_1(T,x)-u_2(T,x)\leq
\Phi(x)-\Phi(x)=0,$ for all $x\in {\mathbb{R}}^n.$\ Let now
$(t_0,x_0)\in [0,T) \times{\mathbb{R}}^n$\ and $\varphi\in
C^{3}([0,T]\times{\mathbb{R}}^n)$\ be such that $w-\varphi$
achieves a strict global maximum at $(t_0,x_0)$. For proving the
theorem it suffices to show that \be\begin{array}{lll}
     &\!\!\!\!\!{\rm min}\{\omega(t_0,x_0),-\frac{\partial }{\partial t} \varphi(t_0,x_0) - \sup_{u \in
U, v \in V}( \frac{1}{2}tr(\sigma\sigma^{T}(t_0, x_0,
 u, v)D^2\varphi(t_0,x_0)\\
 &\!\!\!\!\!\mbox{ }\hskip1cm  + D\varphi.(t_0,x_0)b(t_0,x_0, x, u, v)+ K|\varphi(t_0,x_0)|+K|D\varphi (t_0,x_0) .\sigma(t_0, x_0, u, v)|) \}\le 0.
\end{array}\ee

For this end, applying the method of the separation of variables we introduce the function
$$\Phi_{\varepsilon,\alpha}(t,x,s,y)=u_1(t,x)-u_2(s,y)-\frac{|x-y|^2}{\varepsilon^2}
-\frac{(t-s)^2}{\alpha^2}-\varphi(t,x),$$
where $\varepsilon$ and $\alpha$\ are positive parameters which are devoted to tend to zero.

Since $(t_0, x_0)$\ is a strict global maximum point of $w-\varphi$,
there exists a sequence $(\bar{t},\bar{x},\bar{s},\bar{y})$\ such
that \medskip

(i) $(\bar{t},\bar{x},\bar{s},\bar{y})$\ is a global maximum point
of $\Phi_{\varepsilon, \alpha}$\ in
 $[0,T]\times\bar{B}_r\times\bar{B}_r$\ where $B_r$ is a ball with a
large radius $r$;

\medskip
(ii) $(\bar{t},\bar{x}),\ (\bar{s},\bar{y})\to (t_0,x_0)$ as
$(\varepsilon, \alpha)\to 0$;

\medskip
(iii) $\frac{|\bar{x}-\bar{y}|^2}{\varepsilon^2},\
\frac{(\bar{t}-\bar{s})^2}{\alpha^2}$\ are bounded and tend to zero
when $(\varepsilon, \alpha)\to 0$.\medskip

\noindent Since $u_2$\ is lower semicontinuous we have
$\liminf_{(\varepsilon, \alpha)\rightarrow
0}u_2(\bar{s},\bar{y})\geq u_2(t_0, x_0)$; and thanks to the upper semicontinuity of $u_1$ we have $\limsup_{(\varepsilon, \alpha)\rightarrow
0}u_1(\bar{t},\bar{x})\leq u_1(t_0, x_0)$. On the other hand, from
$\Phi_{\varepsilon,\alpha}(\bar{t},\bar{x},\bar{s},\bar{y})\geq
\Phi_{\varepsilon,\alpha}(t_0,x_0,t_0,x_0)$\ we get
$$u_2(\bar{s},\bar{y})\leq u_1(\bar{t},\bar{x})-u_1(t_0,x_0)+u_2(t_0,x_0)+\varphi(t_0,x_0)-\varphi(\bar{t},\bar{x})
-\frac{|\bar{x}-\bar{y}|^2}{\varepsilon^2}-\frac{(\bar{t}-\bar{s})^2}{\alpha^2}.$$
This yields $\limsup_{(\varepsilon, \alpha)\rightarrow
0}u_2(\bar{s},\bar{y})\leq u_2(t_0, x_0)$. Therefore, we have

\medskip
(iv) $\lim_{(\varepsilon, \alpha)\rightarrow 0}u_2(\bar{s},\bar{y})=
u_2(t_0, x_0)$.
\medskip

\noindent Analogously, we also get

\medskip
(v) $\lim_{(\varepsilon, \alpha)\rightarrow 0}u_1(\bar{t},\bar{x})=
u_1(t_0, x_0)$.
\medskip

\noindent Since $(\bar{t}, \bar{x}, \bar{s}, \bar{y})$\ is a local
maximum point of $\Phi_{\varepsilon, \alpha}$, $u_2({s},
{y})+\frac{|\bar{x}-y|^2}{\varepsilon^2}+\frac{(\bar{t}-s)^2}{\alpha^2}
$\ achieves in $(\bar{s}, \bar{y})$\ a local minimum and from the
definition of a viscosity supersolution of equation (4.1) we have\
$u_2(\bar{s},\bar{y})\geq h(\bar{s},\bar{y}).$ From (iv) we get
$u_2(t_0, x_0)\geq h(t_0, x_0).$\ If $\omega(t_0,x_0)\le 0$,
relation (5.3) is trivially fulfilled. So let us suppose that
$\omega(t_0,x_0)>0$. In this case, we have $h(t_0,x_0)\le
u_2(t_0,x_0)<u_1(t_0,x_0)$. Then according to (v), and since $h$ is
continuous we have

\medskip
(vi) $u_1(\bar{t}, \bar{x})>h(\bar{t}, \bar{x}),$\ for
$\varepsilon>0$\ and $\alpha>0$\ sufficiently small.

\noindent Similarly, since $u_1({t},
{x})-\frac{|x-\bar{y}|^2}{\varepsilon^2}-\frac{(t-\bar{s})^2}{\alpha^2}-\varphi({t},{x})
$\ achieves in $(\bar{t}, \bar{x})$\ a local maximum and from the
definition of a viscosity subsolution of equation (4.1) and (vi) we
have $u_1(\bar{t}, \bar{x})\leq h'(\bar{t}, \bar{x}).$\ From (v) we
get $u_1(t_0,x_0)\leq h'(t_0,x_0).$ Therefore now we get
$u_2(t_0,x_0)< h'(t_0,x_0).$ Similarly, according to (iv) and $h'$
is continuous we have

\medskip
(vii) $u_2(\bar{s}, \bar{y})<h'(\bar{s}, \bar{y}),$\ for
$\varepsilon>0$\ and $\alpha>0$\ sufficiently small.

\medskip

The properties (i) to (vii) and the fact that $u_1 $\ is a viscosity
subsolution and $u_2$\ a viscosity supersolution of equation (5.1)
allow to proceed in the rest of the proof of this lemma exactly as
in the proof of Lemma 3.7 in~\cite{BBE} (our situation here is even
simpler because, contrary to Lemma 3.7 in~\cite{BBE}, we don't have
any integral part in equation (5.1)). So we get:
$$\begin{array}{ll}
 -\frac{\partial \varphi}{\partial t}(t_0,x_0)
                                  &-\sup_{u\in U, v\in V}\left\{\frac{1}{2}tr\left((\sigma\sigma^T)(t_0,x_0,u,v)D^2
                                  \varphi(t_0,x_0)\right)\right.+D\varphi(t_0,x_0)b(t_0,x_0,u,v)\\
&\quad  \left.
+K|\omega(t_0,x_0)|+K|D\varphi(t_0,x_0)\sigma(t_0,x_0,u,v)|\right\}\leq
0,
\end{array}$$
from which relation (5.3) follows easily. Therefore $\omega$ is a viscosity subsolution of
equation (5.2) and  the proof is complete.
\endpf

\medskip

\noindent We now can establish the following comparison principle which is the key for the uniqueness for equation (5.1).\\

\bt We assume that (H3.1) and (H3.2) hold. Let an upper
semicontinuous function $u_1$ (resp., a lower semicontinuous
function $u_2$) $\in \Theta$\ be a viscosity subsolution (resp.,
supersolution) of equation (5.1). Then we have
     \be
     u_1 (t,x) \leq u_2 (t,x) , \hskip 0.5cm \mbox{for all}\ \ (t,x) \in [0,T] \times {\mathbb{R}}^n .
    \ee
If, in particular, both $u_1$ and $u_2$ are continuous viscosity solutions from the class $\Theta$ then they coincide on $[0,T] \times {\mathbb{R}}^n$.
\et

\noindent{\textbf{Proof}.} Theorem 5.1 in~\cite{BL2} establishes a
comparison principle for Hamilton-Jacobi-Bellman equations with
obstable of type (5.2). Such an equation is related to controlled
BSDE with one reflecting barrier, studied in that paper. Letting
$\omega_1=u_1-u_2$ we know from Lemma 5.1 that $\omega_1$ is a
viscosity subsolution of equation (5.2). On the other hand,
$\omega_2=0$ is, obviously, a viscosity solution and, hence, also a
viscosity supersolution of equation (5.2). Both functions $\omega_1$
and $\omega_2$ are in $\Theta$, and the comparison principle stated
in Theorem 5.1 of~\cite{BL2} yields that
$u_1-u_2=\omega_1\le\omega_2=0$, i.e., $u_1\le u_2$ on $[0,T]\times
R^n$. Finally, if $u_1,u_2$ are viscosity solutions of (5.2), they
are both viscosity sub- and supersolution, and from the just proved
comparison result we get the equality of $u_1$ and $u_2$.
\endpf

\br Obviously, since
$\widetilde{W}(t,x)=\lim_{m\rightarrow\infty}\downarrow
W_m(t,x)(\geq W(t,x))$,\ and
$\widetilde{\overline{W}}(t,x)=\lim_{m\rightarrow\infty}\uparrow
\overline{W}_m(t,x)(\leq W(t,x))$\ (for their definitions, see Lemma
4.1 and Remarks 4.2 and 4.3), are a viscosity subsolution and a
supersolution, respectively (see Proposition 4.1 and 4.2), and both
are of linear growth, we have due to Theorem 5.1 that
$\widetilde{W}(t,x)=W(t,x)=\widetilde{\overline{W}}(t,x),\,
(t,x)\in[0,T]\times{\mathbb{R}}^n.$\ Consequently, $W$ is a
viscosity solution of (5.1), and it is unique in the class
$\tilde{\Theta}$. Similarly we get that the upper value function
$U(t,x)$\ is the unique viscosity solution in $\tilde{\Theta}$ of
equation (4.2).

Let us also remark that, since $H^-\leq H^+$ on $[0,T]\times{\mathbb{R}}^n$, any viscosity solution of equation (4.2) is
a supersolution of equation (5.1). Then, again from Theorem
5.1, it follows that $W\leq U$. This justifies calling $W$\ lower
value function and $U$\ upper value function.\er

\br If the Isaacs' condition holds, that is, if for all $(t, x, y,
p, X)\in [0, T]\times {\mathbb{R}}^n \times {\mathbb{R}}\times
{\mathbb{R}}^n\times {\mathbb{S}}^n ,$
$$H^-(t, x, y, p, X)=H^+(t, x, y, p, X),$$
then the equations (5.1) and (4.2) coincide, and from the uniqueness
of the viscosity solution in $\Theta$\ it follows that the lower
value function $W(t,x)$ equals to the upper value function
$U(t,x),$\ that means, the associated stochastic differential game
with reflection has a value.\er

\section{\large{Appendix}}

\vskip0.3cm
 \noindent{\bf\it \textbf{6.1 RBSDES with two Barriers Associated with Forward SDEs}} \vskip0.3cm

 \hskip1cm In this section we give an overview over basic results on RBSDEs with two barriers associated
 with Forward SDEs (for short: FSDEs). We consider measurable functions $b:[0,T]\times \Omega\times
{\mathbb{R}}^n\rightarrow {\mathbb{R}}^n \ $ and
         $\sigma:[0,T]\times \Omega\times {\mathbb{R}}^n\rightarrow {\mathbb{R}}^{n\times d}$
which are supposed to satisfy the following conditions:
 $$
  \begin{array}{ll}
\mbox{(i)}&b(\cdot,x)\ \mbox{and}\ \sigma(\cdot,x)\ \mbox{are} \
{\mathbb{F}}\mbox{-adapted processes, and there exists some}\\
 & \mbox{constant}\ C>0\  \mbox{such that}\\
 &\hskip 1cm|b(t,x)|+|\sigma(t,x)|\leq C(1+|x|), a.s.,\
                                  \mbox{for all}\ 0\leq t\leq T,\ x\in {\mathbb{R}}^n;\\
\mbox{(ii)}&b\ \mbox{and}\ \sigma\ \mbox{are Lipschitz in}\ x,\ \mbox{i.e., there is some constant}\ C>0\ \mbox{such that}\\
           &\hskip 1cm|b(t,x)-b(t,x')|+|\sigma(t,x)-\sigma(t,x')|\leq C| x-x'|,\ a.s.,\\
 & \hbox{ \ \ }\hskip7cm\mbox{for all}\ 0\leq t \leq T,\ x,\ x'\in {\mathbb{R}}^n.\\
 \end{array}
  \eqno{\mbox{(H6.1)}}
  $$\par
  We now consider the following SDE parameterized by the
  initial condition $(t,\zeta)\in[0,T]\times L^2(\Omega,{\cal{F}}_t,P;{\mathbb{R}}^n)$:
  \be
  \left\{
  \begin{array}{rcl}
  dX_s^{t,\zeta}&=&b(s,X_s^{t,\zeta})ds+\sigma(s,X_s^{t,\zeta})dB_s,\ s\in[t,T],\\
  X_t^{t,\zeta}&=&\zeta.
  \end{array}
  \right.
  \ee
Under the assumption (H6.1), SDE (6.1) has a unique strong solution
and, for any $p\geq 2,$\ there exists $C_{p}\in {\mathbb{R}}$\ such
that, for any $t\in[0,T]\ \mbox{and}\ \zeta,\zeta'\in
L^p(\Omega,{\cal{F}}_t,P;{\mathbb{R}}^n),$
 \be
 \begin{array}{rcl}
 E[\sup\limits_{t\leq s\leq T}| X_s^{t,\zeta}-X_s^{t,\zeta'}|^p|{\cal{F}}_t]
                             &\leq& C_{p}|\zeta-\zeta'|^p, \ \ a.s.,\\
  E[\sup\limits_{t\leq s\leq T}| X_s^{t,\zeta}|^p|{\cal{F}}_t]
                       &\leq& C_{p}(1+|\zeta|^p),\ \  a.s.
 \end{array}
\ee \noindent These well-known standard estimates can be consulted,
for instance, in Ikeda, Watanabe~\cite{IW}, pp.166-168, and also in
Karatzas, Shreve~\cite{KSH}, pp.289-290. We emphasize that the
constant $C_{p}$ in (6.2) only depends on the Lipschitz and the
growth constants of $b$ and $\sigma$.

 Let now be given three real valued functions $f(t,x,y,z)$, $\Phi(x)$\ and $h(t,x), h'(t,x)$\ which shall satisfy the
following conditions:
$$
\begin{array}{ll}
\mbox{(i)}&\Phi:\Omega\times {\mathbb{R}}^n\rightarrow {\mathbb{R}}
\ \mbox{is an}\ {\cal{F}}_T\otimes{\cal{B}}({\mathbb{R}}^n)
             \mbox{-measurable random variable and}\\
          &\hskip 0.5cm f:[0,T]\times \Omega\times {\mathbb{R}}^n\times {\mathbb{R}}\times
          {\mathbb{R}}^d \rightarrow {\mathbb{R}},\ \
h, h':\Omega\times [0, T]\times {\mathbb{R}}^n\rightarrow {\mathbb{R}}\ \mbox{}\\
          & \mbox{are measurable processes such that, }\\
          &\hskip 0.5cm f(\cdot,x,y,z),\ h(\cdot,x),\ h'(\cdot,x)\ \mbox{are}\
          {\mathbb{F}}\mbox{-adapted, for all $(x, y, z)\in{\mathbb{R}}^n\times {\mathbb{R}}\times
          {\mathbb{R}}^d $;}\\
\end{array}$$
$$\begin{array}{ll}
 \mbox{(ii)}&\mbox{There exist constants}\ \mu>0\ \mbox{such that, P-a.s.,}\\
          &| f(t,x,y,z)-f(t,x',y',z')|\leq \mu(|x-x'|+ |y-y'|+|z-z'|);\\
          &| \Phi(x)-\Phi(x')|\leq \mu|y-y'|;\\
          &|h(t,x)-h(t,x')|+ |h'(t,x)-h'(t,x')|\leq \mu|x-x'|;\\
&\hskip 3cm \mbox{for all}\ 0\leq t\leq T,\ x,\ x'\in
{\mathbb{R}}^n,\ y,\ y'\in {\mathbb{R}}\ \mbox{and}\ z,\ z'\in
{\mathbb{R}}^d;\\
 \mbox{(iii)}&f\ \mbox{and}\ \Phi \ \mbox{satisfy a linear growth condition, i.e., there exists some}\ C>0\\\
    & \mbox{such that, dt}\times \mbox{dP-a.e.},\ \mbox{for all}\ x\in
    {\mathbb{R}}^n,\\
    &\hskip 2cm|f(t,x,0,0)| + |\Phi(x)|\leq C(1+|x|)\\
&\mbox{and, moreover,}\\
&\hskip 0.5cm h(\cdot,x),\ h'(\cdot,x)\ \mbox{are continuous in}\ t,
\ h(t, x)< h'(t, x),\ \mbox{a.s.}, \mbox{and}\\
 &\hskip 0.5cm h(T, x)\leq \Phi(x)\leq h'(T, x),\
a.s.,\ \mbox{for all}\ (t, x)\in [0, T]\times{\mathbb{R}}^n.\\
\end{array}
\eqno{\mbox{(H6.2)}} $$

 With the help of the above assumptions we can verify that the coefficient $f(s,X_s^{t,\zeta},y,z)$\ satisfies the hypotheses (A1),
 (A2),  $\xi:=\Phi(X_T^{t,\zeta})$ $\in
 L^2(\Omega,{\cal{F}}_T,P;{\mathbb{R}})$\ and $L_s:=h(s,X_s^{t,\zeta}),\ U_s:=h'(s,X_s^{t,\zeta})\in {\cal{A}}_{c}^{2}(0,T; {\mathbb{R}})$.
 Therefore, according to Lemma 2.6 the following RBSDE with two
 barriers possesses a unique solution:
\be
\begin{array}{lll}
 &{\rm (i)}Y^{t,\zeta} \in {\cal{S}}^2(0,
T; {\mathbb{R}}),\ Z^{t,\zeta} \in
{\cal{H}}^{2}(0,T;{\mathbb{R}}^{d})\  \mbox{and}\ \
  K^{+,t,\zeta},\ K^{-,t,\zeta} \in {\cal{A}}^2_{c}(0,
T; {\mathbb{R}});\\
&{\rm (ii)} Y^{t,\zeta}_s = \Phi(X_T^{t,\zeta}) +
\int_s^Tf(r,X_r^{t,\zeta},Y^{t,\zeta}_r,Z^{t,\zeta}_r)dr +
(K^{+,t,\zeta}_{T} -K^{+,t,\zeta}_{s})-(K^{-,t,\zeta}_{T} -K^{-,t,\zeta}_{s})\\
&\ \hskip1.7cm - \int^T_sZ^{t,\zeta}_rdB_r,\  s\in [t,T];\ \\
&{\rm(iii)}h(s, X_s^{t,\zeta})\leq Y^{t,\zeta}_s \leq h'(s,
X_s^{t,\zeta}),\ \ \mbox{a.s.},\
\mbox{for any}\ s\in [t,T];\\
&{\rm (iv)} \int_t^T(Y^{t,\zeta}_r -
h(r,X_r^{t,\zeta}))dK^{+,t,\zeta}_{r}=\int_t^T(h'(r,X_r^{t,\zeta})-Y^{t,\zeta}_r)dK^{-,t,\zeta}_{r}=0.
\end{array}
\ee

\bp  We suppose that the hypotheses (H6.1) and (H6.2) hold. Then,
for any\, $0\leq t\leq T$ and the associated initial conditions
 $\zeta,\zeta'\in L^2(\Omega,{\cal{F}}_t,P;{\mathbb{R}}^n)$, we
have the following estimates:\\  $\mbox{}\hskip3cm\mbox{\rm(i)}
E[\sup\limits_{t\leq s\leq T}|Y_s^{t,\zeta}|^2 |{{\cal{F}}_t}]\leq
C(1+|\zeta|^2),\
a.s.; $\\
$\mbox{}\hskip3cm\mbox{\rm(ii)}E[\sup\limits_{t\leq s\leq
T}|Y_s^{t,\zeta}-Y_s^{t,\zeta'}|^2|{{\cal{F}}_t}] \leq C|\zeta-\zeta'|^2,\  a.s. $\\
In particular, \be
 \begin{array}{lll}
\mbox{\rm(iii)}&|Y_t^{t,\zeta}|\leq C(1+|\zeta|),\  a.s.; \hskip3cm\\
\mbox{\rm(iv)}&|Y_t^{t,\zeta}-Y_t^{t,\zeta'}|\leq C|\zeta-\zeta'|,\  a.s.\\
\end{array}
\ee The above constant $C>0$\ depends only on the Lipschitz and the
growth constants of $b$,\ $\sigma$, $f$, $\Phi$\ and $h$. \ep
 \noindent \textbf{Proof}. From Lemma 2.8 combined with (6.2) we obtain easily (i). So we need only to
 prove (ii). For an arbitrarily fixed $\varepsilon>0$, we
 define the function $\psi_\varepsilon(x)=(|x|^2+\varepsilon)^{\frac{1}{2}},\ x\in
 {\mathbb{R}}^n.$\ Obviously, $|x|\leq \psi_\varepsilon(x)\leq |x|+\varepsilon^{\frac{1}{2}},\ x\in
 {\mathbb{R}}^n.$\ Furthermore, for all $x\in {\mathbb{R}}^n,$
 $$
D\psi_\varepsilon(x)=\frac{x}{(|x|^2+\varepsilon)^{\frac{1}{2}}},\ \
\ \ \
D^2\psi_\varepsilon(x)=\frac{I}{(|x|^2+\varepsilon)^{\frac{1}{2}}}-\frac{x\otimes
x}{(|x|^2+\varepsilon)^{\frac{3}{2}}}.
 $$
Then, we have \be |D\psi_\varepsilon(x)|\leq 1,\ \ \
|D^2\psi_\varepsilon(x)||x|\leq
\frac{C}{(|x|^2+\varepsilon)^{\frac{1}{2}}}|x|\leq C,\ \ x\in
{\mathbb{R}}^n,\ee for some constant $C$\ that doesn't dependent on
$\varepsilon$. Let us denote by $X^{t,\zeta}$\ and $X^{t,\zeta'}$\
the unique solution of SDE (6.1) with initial data $(t, \zeta)$\ and
$(t, \zeta')$, respectively. Moreover, recall that $\mu$\ is the
Lipschitz constant of $h,\ h', \ \Phi,\ \ f$. We consider the
following two RBSDEs: \be
\begin{array}{lll}
 &{\rm (i)}\,  \widetilde{Y}\in
{\cal{S}}^2(0, T; {\mathbb{R}}), \ \widetilde{Z}\in
{\cal{H}}^{2}(0,T;{\mathbb{R}}^{d})\ \ \mbox{and}\ \
  \widetilde{K}^+,\ \widetilde{K}^-\in {\cal{A}}^2_{c}(0, T; {\mathbb{R}});\\
&{\rm (ii)}\,  \widetilde{Y}_s = \Phi(X_T^{t,\zeta}) +
\mu\psi_\varepsilon(X_T^{t,\zeta}-X_T^{t,\zeta'})+
\int_s^T(f(r,X_r^{t,\zeta},\widetilde{Y}_r,\widetilde{Z}_r)+\mu|X_r^{t,\zeta}-X_r^{t,\zeta'}|)dr\\
&\ \hskip3cm + (\widetilde{K}^+_{T} -
\widetilde{K}^+_{s})-(\widetilde{K}^-_{T}
-\widetilde{K}^-_{s}) - \int^T_s\widetilde{Z}_rdB_r,\ s\in [t,T];\ \\
&{\rm(iii)}\,
h(s,X_s^{t,\zeta})+\mu\psi_\varepsilon(X_s^{t,\zeta}-X_s^{t,\zeta'})\leq\widetilde{Y}_s
\leq
h'(s,X_s^{t,\zeta})+\mu\psi_\varepsilon(X_s^{t,\zeta}-X_s^{t,\zeta'}), \mbox{a.s.}, s\in [t,T];\\
&{\rm (iv)}\,  \int_t^T(\widetilde{Y}_r -
h(r,X_r^{t,\zeta})-\mu\psi_\varepsilon(X_r^{t,\zeta}-X_r^{t,\zeta'}))d\widetilde{K}^+_{r}
=\int_t^T(h'(r,X_r^{t,\zeta})+\mu\psi_\varepsilon(X_r^{t,\zeta}-X_r^{t,\zeta'})\\
&\ \hskip0.5cm -\widetilde{Y}_r)d
\widetilde{K}^{-}_{r}=0.\end{array} \ee and \be
\begin{array}{lll}
 &{\rm (i)}\,  \bar{Y}\in
{\cal{S}}^2(0, T; {\mathbb{R}}), \ \bar{Z} \in
{\cal{H}}^{2}(0,T;{\mathbb{R}}^{d})\ \ \mbox{and}\ \
  \bar{K}^+,\ \bar{K}^-\in {\cal{A}}^2_{c}(0, T; {\mathbb{R}});\\
&{\rm (ii)}\, \bar{Y}_s = \Phi(X_T^{t,\zeta}) -
\mu|X_T^{t,\zeta}-X_T^{t,\zeta'}|+
\int_s^T(f(r,X_r^{t,\zeta},\bar{Y}_r,\bar{Z}_r)-\mu|X_r^{t,\zeta}-X_r^{t,\zeta'}|)dr\\
&\ \hskip3cm + (\bar{K}^+_{T} - \bar{K}^+_{s}) - (\bar{K}^-_{T} -
\bar{K}^-_{s}) -\int^T_s\bar{Z}_rdB_r,\ \ \ \  s\in [t,T];\ \\
&{\rm(iii)}\,
h(s,X_s^{t,\zeta})-\mu\psi_\varepsilon(X_s^{t,\zeta}-X_s^{t,\zeta'})\leq\bar{Y}_s
\leq
h'(s,X_s^{t,\zeta})-\mu\psi_\varepsilon(X_s^{t,\zeta}-X_s^{t,\zeta'}), \mbox{a.s.}, s\in [t,T];\\
&{\rm (iv)}\, \int_t^T(\bar{Y}_r -
h(r,X_r^{t,\zeta})+\mu\psi_\varepsilon(X_r^{t,\zeta}-X_r^{t,\zeta'}))d\bar{K}^+_{r}=\int_t^T(
h'(r,X_r^{t,\zeta})-\mu\psi_\varepsilon(X_r^{t,\zeta}-X_r^{t,\zeta'})\\
& \ \hskip0.5cm -\bar{Y}_r)d\bar{K}^-_{r}=0.
\end{array}\ee

\noindent Obviously, their coefficients satisfy the assumptions in
(H6.2) and they admit unique solutions $(\widetilde{Y},
\widetilde{Z}, \widetilde{K}^+, \widetilde{K}^-)$\ and $(\bar{Y},
\bar{Z}, \bar{K}^+,\bar{K}^-)$, respectively. Moreover, from the
comparison theorem for RBSDEs with two barriers (Lemma 2.7)
\be\bar{Y}_s\leq Y_s^{t,\zeta}\leq \widetilde{Y}_s,\ \ \
\bar{Y}_s\leq Y_s^{t,\zeta'}\leq \widetilde{Y}_s,\ \ \mbox{P-a.s.,
for all}\ s\in [t, T].\ee We shall still introduce two other RBSDEs
with two reflecting barriers: \be
\begin{array}{lll}
 &{\rm (i)}\,  \widetilde{Y}'\in
{\cal{S}}^2(0, T; {\mathbb{R}}), \ \widetilde{Z}'\in
{\cal{H}}^{2}(0,T;{\mathbb{R}}^{d})\  \mbox{and}\ \
  \widetilde{K}^{+'},\ \widetilde{K}^{-'}\in{\cal{A}}^2_{c}(0, T; {\mathbb{R}});\\
&{\rm (ii)}\, \widetilde{Y}'_s = \Phi(X_T^{t,\zeta}) +\\
&\ \ \ \
\int_s^T[f(r,X_r^{t,\zeta},\widetilde{Y}'_r+\mu\psi_\varepsilon(X_r^{t,\zeta}-X_r^{t,\zeta'}),\widetilde{Z}'_r
+\mu D\psi_\varepsilon(X_r^{t,\zeta}-X_r^{t,\zeta'})(\sigma(r,X_r^{t,\zeta})-\sigma(r,X_r^{t,\zeta'})))\\
&\ \hskip0.5cm +\mu|X_r^{t,\zeta}-X_r^{t,\zeta'}|+\mu D\psi_\varepsilon(X_r^{t,\zeta}-X_r^{t,\zeta'})(b(r,X_r^{t,\zeta})-b(r,X_r^{t,\zeta'}))\\
&\hskip1.5cm+\frac{1}{2}\mu(D^2\psi_\varepsilon(X_r^{t,\zeta}-X_r^{t,\zeta'})
(\sigma(r,X_r^{t,\zeta})-\sigma(r,X_r^{t,\zeta'})),\sigma(r,X_r^{t,\zeta})-\sigma(r,X_r^{t,\zeta'}))]dr\\
&\ \ \ \ + (\widetilde{K}^{+'}_{T} - \widetilde{K}^{+'}_{s})
-(\widetilde{K}^{-'}_{T} -\widetilde{K}^{-'}_{s}) - \int^T_s\widetilde{Z}'_rdB_r,\ \ \ \  s\in [t,T];\ \\
&{\rm(iii)}\, h(s, X_s^{t,\zeta})\leq\widetilde{Y}'_s \leq h'(s,
X_s^{t,\zeta}),\ \ \mbox{a.s.},\
 s\in [t,T]; \\
&{\rm (iv)}\, \int_t^T(\widetilde{Y}'_r -
h(r,X_r^{t,\zeta}))d\widetilde{K}^{+'}_{r}=\int_t^T(h'(r,X_r^{t,\zeta})-\widetilde{Y}'_r
)d\widetilde{K}^{-'}_{r}=0,
\end{array} \ee
\noindent and \be
\begin{array}{lll}
 &{\rm (i)}\,  \bar{Y}'\in
{\cal{S}}^2(0, T; {\mathbb{R}}), \ \bar{Z}'\in
{\cal{H}}^{2}(0,T;{\mathbb{R}}^{d})\  \mbox{and}\ \
  \bar{K}^{+'},\ \bar{K}^{-'}\in{\cal{A}}^2_{c}(0, T; {\mathbb{R}});\\
&{\rm (ii)}\, \bar{Y}'_s = \Phi(X_T^{t,\zeta})-
\mu|X_T^{t,\zeta}-X_T^{t,\zeta'}|+\mu\psi_\varepsilon(X_T^{t,\zeta}-X_T^{t,\zeta'})+\\
&\int_s^T[f(r,X_r^{t,\zeta},\bar{Y}'_r-\mu\psi_\varepsilon(X_r^{t,\zeta}-X_r^{t,\zeta'}),\bar{Z}'_r
-\mu D\psi_\varepsilon(X_r^{t,\zeta}-X_r^{t,\zeta'})(\sigma(r,X_r^{t,\zeta})-\sigma(r,X_r^{t,\zeta'})))\\
&\ \hskip0.5cm -\mu|X_r^{t,\zeta}-X_r^{t,\zeta'}|-\mu D\psi_\varepsilon(X_r^{t,\zeta}-X_r^{t,\zeta'})(b(r,X_r^{t,\zeta})-b(r,X_r^{t,\zeta'}))\\
&\
\hskip1.5cm-\frac{1}{2}\mu(D^2\psi_\varepsilon(X_r^{t,\zeta}-X_r^{t,\zeta'})
(\sigma(r,X_r^{t,\zeta})-\sigma(r,X_r^{t,\zeta'})),\sigma(r,X_r^{t,\zeta})-\sigma(r,X_r^{t,\zeta'}))]dr\\
&\ \hskip0.5cm + (\bar{K}^{+'}_{T} - \bar{K}^{+'}_{s}) -
(\bar{K}^{-'}_{T} -\bar{K}^{-'}_{s}) -\int^T_s\bar{Z}'_rdB_r,\ \ \ \  s\in [t,T];\ \\
&{\rm(iii)}h(s, X_s^{t,\zeta})\leq\bar{Y}'_s \leq h'(s,
X_s^{t,\zeta}),\ \ \mbox{a.s.},\ \ s\in [t,T];\\
&{\rm (iv)} \int_t^T(\bar{Y}'_r -
h(r,X_r^{t,\zeta}))d\bar{K}^{+'}_{r}= \int_t^T(
h'(r,X_r^{t,\zeta})-\bar{Y}'_r)d\bar{K}^{-'}_{r}=0.\\
\end{array} \ee

Obviously, also the RBSDEs with two barriers (6.9) and (6.10)
satisfy the assumption (H6.2) and, thus, admit unique solutions
$(\widetilde{Y}', \widetilde{Z}', \widetilde{K}^{+'},
\widetilde{K}^{-'})$\ and $(\bar{Y}', \bar{Z}', \bar{K}^{+'},
\bar{K}^{-'})$, respectively. On the other hand, from the uniqueness
of the solution of RBSDE with two barriers we know that
\be\begin{array}{lll}
&\widetilde{Y}'_s=\widetilde{Y}_s-\mu\psi_\varepsilon(X_s^{t,\zeta}-X_s^{t,\zeta'}),\
\mbox{for all}\ s\in [t, T],\ \mbox{P-a.s.,}\\
&\widetilde{Z}'_s=\widetilde{Z}_s-\mu
D\psi_\varepsilon(X_s^{t,\zeta}-X_s^{t,\zeta'})(\sigma(s,X_s^{t,\zeta}
)-\sigma(s,X_s^{t,\zeta'})),\  \mbox{dsdP-a.e. on}\ [t,
T]\times\Omega,\\
&\widetilde{K}^{+'}_s=\widetilde{K}^{+}_s,\
\widetilde{K}^{-'}_s=\widetilde{K}^{-}_s\ \mbox{for all}\ s\in [t,
T],\ \mbox{P-a.s.}\end{array}\ee \noindent and \be\begin{array}{lll}
&\bar{Y}'_s=\bar{Y}_s+\mu\psi_\varepsilon(X_s^{t,\zeta}-X_s^{t,\zeta'}),\
\mbox{for all}\ s\in [t, T],\ \mbox{P-a.s.,}\\
&\bar{Z}'_s=\bar{Z}_s+\mu
D\psi_\varepsilon(X_s^{t,\zeta}-X_s^{t,\zeta'})(\sigma(s,X_s^{t,\zeta}
)-\sigma(s,X_s^{t,\zeta'})),\  \mbox{dsdP-a.e. on}\ [t,
T]\times\Omega,\\
&\bar{K}^{+'}_s=\bar{K}^{+}_s,\ \bar{K}^{-'}_s=\bar{K}^{-}_s,
\mbox{for all}\ s\in [t, T],\ \mbox{P-a.s.}\end{array}\ee

\noindent Then, for the notations introduced in Lemma 2.9 we have
\be\begin{array}{lll} &\Delta g(r, \widetilde{Y}'_r,
\widetilde{Z}'_r)\\
&=f(r,X_r^{t,\zeta},\widetilde{Y}'_r+\mu\psi_\varepsilon(X_r^{t,\zeta}-X_r^{t,\zeta'}),\widetilde{Z}'_r
+\mu D\psi_\varepsilon(X_r^{t,\zeta}-X_r^{t,\zeta'})(\sigma(r,X_r^{t,\zeta})-\sigma(r,X_r^{t,\zeta'})))\\
&-f(r,X_r^{t,\zeta},\widetilde{Y}'_r-\mu\psi_\varepsilon(X_r^{t,\zeta}-X_r^{t,\zeta'}),\widetilde{Z}'_r
-\mu D\psi_\varepsilon(X_r^{t,\zeta}-X_r^{t,\zeta'})(\sigma(r,X_r^{t,\zeta})-\sigma(r,X_r^{t,\zeta'})))\\
&\ \hskip0.5cm +2\mu|X_r^{t,\zeta}-X_r^{t,\zeta'}|+2\mu D\psi_\varepsilon(X_r^{t,\zeta}-X_r^{t,\zeta'})(b(r,X_r^{t,\zeta})-b(r,X_r^{t,\zeta'}))\\
&\ \hskip1.5cm+\mu(D^2\psi_\varepsilon(X_r^{t,\zeta}-X_r^{t,\zeta'})
(\sigma(r,X_r^{t,\zeta})-\sigma(r,X_r^{t,\zeta'})),\sigma(r,X_r^{t,\zeta})-\sigma(r,X_r^{t,\zeta'}));\\
&\Delta\xi=\mu|X_T^{t,\zeta}-X_T^{t,\zeta'}|-\mu\psi_\varepsilon(X_T^{t,\zeta}-X_T^{t,\zeta'}).\\
\end{array}\ee
From (6.5) and the Lipschitz continuity of $f,\  b$\ and $\sigma$\
we get
$$\begin{array}{lll}
&|\Delta g(r, \widetilde{Y}'_r, \widetilde{Z}'_r)|\leq
C|X_r^{t,\zeta}-X_r^{t,\zeta'}|+C\varepsilon^{\frac{1}{2}},\
\mbox{P-a.s.},\\
&|\Delta \xi|\leq
C|X_T^{t,\zeta}-X_T^{t,\zeta'}|+C\varepsilon^{\frac{1}{2}},\
\mbox{P-a.s.},\end{array}$$\ where the constant $C$ \ is independent
of $\varepsilon$. Therefore, from Lemma 2.9 and (6.2) we get
$$E[\sup_{t\leq s\leq T}|\widetilde{Y}'_s-\bar{Y}'_s|^2|{\cal F}_t]
\leq C|\zeta-\zeta'|^2+C\varepsilon,\ \mbox{P-a.s.}$$
Furthermore, from (6.8), (6.11), (6.12) and (6.2) we have
$$\begin{array}{rcl}
& & E[\sup_{t\leq s\leq T}|Y^{t,\zeta}_s-Y^{t,\zeta'}_s|^2|{\cal
F}_t]\leq E[\sup_{t\leq s\leq T}|\widetilde{Y}_s-\bar{Y}_s|^2|{\cal F}_t]\\
& &\leq 2E[\sup_{t\leq s\leq T}|\widetilde{Y}'_s-\bar{Y}'_s|^2|{\cal
F}_t]+16\mu^2(E[\sup_{t\leq s\leq
T}|X^{t,\zeta}_s-X^{t,\zeta'}_s|^2|{\cal F}_t]+\varepsilon)\\
& &\leq  C|\zeta-\zeta'|^2+C\varepsilon,\ \mbox{P-a.s.}\end{array}$$

\noindent Finally, we let $\varepsilon$\ tend to 0, and we get (ii). The
proof is complete.\endpf

 \vskip 0.5cm Let us now introduce the random field:
\be u(t,x)=Y_s^{t,x}|_{s=t},\ (t, x)\in [0, T]\times{\mathbb{R}}^n,
\ee where $Y^{t,x}$ is the solution of the RBSDE with two barriers (6.3),
with $\zeta\in L^2(\Omega,{\cal{F}}_t,P;{\mathbb{R}}^n)$\ being
replaced by $x \in
{\mathbb{R}}^n$.\\
As a consequence of Proposition 6.1 we have that, for all $t \in [0,
T] $, P-a.s.,
 \be
\begin{array}{ll}
\mbox{(i)}&| u(t,x)-u(t,y)| \leq C|x-y|,\ \mbox{for all}\ x, y\in {\mathbb{R}}^n;\\
\mbox{(ii)}&| u(t,x)|\leq C(1+|x|),\ \mbox{for all}\ x\in {\mathbb{R}}^n.\\
\end{array}
\ee

 The random field $u$\ and $Y^{t,\zeta},\ (t, \zeta)\in [0,
T]\times L^2(\Omega,{\cal{F}}_t,P;{\mathbb{R}}^n),$\ are related by
the following theorem.
 \bp Under the assumptions (H6.1) and (H6.2) we have, for any $t\in [0, T]$\ and $\zeta\in
L^2(\Omega,{\cal{F}}_t,P;{\mathbb{R}}^n),$ \be
u(t,\zeta)=Y_t^{t,\zeta},\ \mbox{ P-a.s.} \ee \ep The proof of
Proposition 6.2 can be got by adapting the corresponding argument of
Peng~\cite{Pe1} to RBSDEs, we give it for the reader's convenience.
It makes use of the following definition.

 \noindent
\bde For any t $\in [0, T]$, a sequence $\{A_i\}_{i=1}^{N}\subset
{\cal{F}}_t\ (\mbox{with}\ 1\leq N\leq \infty)$ is called a
partition of $(\Omega, {\cal{F}}_t)$\ if\ \
$\cup_{i=1}^{N}A_i=\Omega$\ and $ A_i\cap A_j=\phi, \
\mbox{whenever}\ i\neq j.$ \ede
 \noindent \textbf{Proof} (of Proposition 6.2): We first consider
the case that $\zeta$ is a simple random variable of the form\
\be\zeta=\sum\limits^N\limits_{i=1}x_i\textbf{1}_{A_i},
                        \ee
where $\{A_i\}^N_{i=1}$\ is a finite partition of $(\Omega,{\cal{F}}_t)$\
and $x_i\in {\mathbb{R}}^n$,\ for $1\leq i\leq N.$\\
For each $i$, we put $(X_s^i,Y_s^i,Z_s^i)\equiv
                 (X_s^{t,x_i},Y_s^{t,x_i},Z_s^{t,x_i}).$ Then $X^i$ is the solution of the SDE
$$
X^i_s =x_i +\int^s_t b(r,X^i_r)dr +\int^s_t \sigma (r,X^i_r)dB_r,\
s\in [t,T],
$$
 and $(Y^i,Z^i, K^{+,i}, K^{-,i})$ is the solution of the associated RBSDE
$$\begin{array}{ll}
&Y^i_s =\Phi(X^i_T) +\int^T_s
f(r,X^i_r,Y^i_r,Z^i_r)dr+(K^{+,i}_T-K^{+,i}_s)-(K^{-,i}_T-K^{-,i}_s)
   -\int^T_s Z^i_r dB_r,\ s\in [t,T],\\
&h(s, X_s^i)\leq Y^i_s \leq h'(s, X_s^i),\ \  \int_t^T(Y^i_r -
h(r,X_r^i))dK^{+,i}_{r}=\int_t^T(h'(r,X_r^i)-Y^i_r)dK^{-,i}_{r}=0.\\
\end{array}
$$
The above two equations are multiplied by $\textbf{1}_{A_i}$\ and
summed up with respect to $i$. Thus, taking into account that
$\sum\limits_i \varphi (x_i)\textbf{1}_{A_i}=\varphi (\sum\limits_i
x_i \textbf{1}_{A_i})$, for any function $\varphi$, we get
$$
\begin{array}{rcl}
 \sum\limits_{i=1}\limits^{N} \textbf{1}_{A_i} X^i_s &=&\sum\limits _{i=1}\limits^{N}
  x_i \textbf{1}_{A_i}+ \int^s_t b(r,\sum\limits _{i=1}\limits^{N} \textbf{1}_{A_i} X^i_r )dr
  +\int^s_t \sigma (r,\sum\limits_{i=1}\limits^{N} \textbf{1}_{A_i} X^i_r)dB_r
\end{array}
$$and
$$
\begin{array}{lll}
& &\sum\limits _{i=1}\limits^{N}\textbf{1}_{A_i} Y^i_s  =
\Phi(\sum\limits _{i=1}\limits^{N} \textbf{1}_{A_i} X^i_T)+\int^T_s
f(r,\sum \limits_{i=1}^{N} \textbf{1}_{A_i} X^i_r,
  \sum\limits _{i=1}^{N} \textbf{1}_{A_i} Y^i_r,
 \sum\limits _{i=1}^{N} \textbf{1}_{A_i} Z^i_r)dr \\
  & &\ \hskip2cm +(\sum\limits _{i=1}^{N} \textbf{1}_{A_i}K_T^{+,i}-
  \sum\limits _{i=1}^{N} \textbf{1}_{A_i}K_s^{+,i})-(\sum\limits _{i=1}^{N}
  \textbf{1}_{A_i}K_T^{-,i}
  -\sum\limits _{i=1}^{N} \textbf{1}_{A_i}K_s^{-,i})-\int^T_s \sum\limits _{i=1}^{N} \textbf{1}_{A_i} Z^i_r
  dB_r,\\
& & h(s, \sum\limits _{i=1}^{N} \textbf{1}_{A_i}X_s^i) \leq
\sum\limits _{i=1}^{N} \textbf{1}_{A_i}Y^i_s \leq h'(s, \sum\limits
_{i=1}^{N} \textbf{1}_{A_i}X_s^i), \\
& & \int_t^T(\sum\limits _{i=1}^{N} \textbf{1}_{A_i}Y^i_r -
h(r,\sum\limits _{i=1}^{N} \textbf{1}_{A_i}X_r^i))d(\sum\limits
_{i=1}^{N} \textbf{1}_{A_i}K^{+,i}_{r})=\int_t^T(h'(r,\sum\limits
_{i=1}^{N} \textbf{1}_{A_i}X_r^i)-\sum\limits _{i=1}^{N}
\textbf{1}_{A_i}Y^i_r)d(\sum\limits _{i=1}^{N} \textbf{1}_{A_i}K^{-,i}_{r})\\
& &=0.\\
\end{array}
$$
Then the strong uniqueness property of the SDE and the associated
RBSDE with two barriers yields: For $s\in [t, T],$
$$
X^{t,\zeta}_s =\sum \limits_{i=1}^{N} X^i_s \textbf{1}_{A_i},\
(Y^{t,\zeta}_s, Z^{t,\zeta}_s, K^{+,t,\zeta}_s, K^{-,t,\zeta}_s)
=(\sum \limits_{i=1}^{N} \textbf{1}_{A_i} Y^i_s, \sum
\limits_{i=1}^{N} \textbf{1}_{A_i} Z^i_s, \sum \limits_{i=1}^{N}
\textbf{1}_{A_i} K^{+,i}_s, \sum \limits_{i=1}^{N} \textbf{1}_{A_i}
K^{-,i}_s).
$$
Finally, from $u(t,x_i)=Y^i_t,\ 1\leq i\leq N$, we deduce that
$$
Y^{t,\zeta}_t=\sum \limits_{i=1}^{N}
Y^i_t\textbf{1}_{A_i}=\sum\limits_{i=1}^{N}u(t,x_i) \textbf{1}_{A_i}
=u(t,\sum \limits_{i=1}^{N} x_i \textbf{1}_{A_i}) =u(t,\zeta).
$$
Therefore, for simple random variables, we have the desired result.

Given a general $\zeta\in L^2 (\Omega ,{\mathcal{F}}_t
,P;{\mathbb{R}}^n)$ we can choose a sequence of simple random
variables $\{\zeta_i\}$ which
 converges to $\zeta$ in $L^2(\Omega ,{\mathcal{F}}_t
,P;{\mathbb{R}}^n)$. Consequently, from the estimates (6.4), (6.15)
and the first step of the proof, we have
$$
\begin{array}{lrcl}
&E|Y^{t,\zeta_i}_t-Y^{t,\zeta}_t|^2&\leq&CE|\zeta_i -\zeta|^2\rightarrow 0,\ i\rightarrow\infty,\\
\mbox{ }\hskip1cm&
E|u(t,\zeta_i)-u(t,\zeta)|^2 &\leq& CE|\zeta_i -\zeta|^2 \rightarrow 0,\ i\rightarrow\infty,\\
\hbox{and}\hskip1cm& Y^{t,\zeta_i}_t&=& u(t,\zeta_i),\ i\geq 1.
\end{array}
$$
Then the proof is complete.\endpf

\vskip0.3cm
 \noindent{\bf\it \textbf{6.2 The proof of Theorem 3.1}} \vskip0.3cm
\noindent{Proof.} To simplify notations we put
$$W_\delta(t,x) =\mbox{essinf}_{\beta \in {\mathcal{B}}_{t,
t+\delta}}\mbox{esssup}_{u \in {\mathcal{U}}_{t,
t+\delta}}G^{t,x;u,\beta(u)}_{t,t+\delta} [W(t+\delta,
X^{t,x;u,\beta(u)}_{t+\delta})].$$ The proof that $W_\delta(t,x)$\
coincides with $W(t,x)$\ will be split into a sequel of lemmata
which all suppose that (H3.1) and (H3.2) are satisfied. Let us fix
$(t, x)\in [0, T]\times {\mathbb{R}}^n.$

\bl $W_\delta(t,x)$\ is deterministic.\el The proof of this lemma
uses the same ideas as that of Proposition 3.1 so that it can be
omitted here.\endpf

\bl$W_\delta(t,x)\leq W(t,x).$\el

\noindent\textbf{ Proof}. Let $\beta\in {\mathcal{B}}_{t, T}$\ be
arbitrarily fixed. Then, given a $u_2(\cdot)\in
{\mathcal{U}}_{t+\delta, T},$\ we define as follows the restriction
$\beta_1$\ of $\beta$\ to ${\mathcal{U}}_{t, t+\delta}:$
$$\beta_1(u_1):=\beta(u_1\oplus u_2 )|_{[t,
t+\delta]},\ \mbox{ }\ u_1(\cdot)\in {\mathcal{U}}_{t, t+\delta},
$$
where $u_1\oplus u_2:=u_1\textbf{1}_{[t,
t+\delta]}+u_2\textbf{1}_{(t+\delta, T]}$\ extends $u_1(\cdot)$\ to
an element of ${\mathcal{U}}_{t, T}$. It is easy to check that
$\beta_1\in {\mathcal{B}}_{t, t+\delta}.$\ Moreover, from the
nonanticipativity property of $\beta$\ we deduce that $\beta_1$\ is
independent of the special choice of $u_2(\cdot)\in
{\mathcal{U}}_{t+\delta, T}.$\ Consequently, from the definition of
$W_\delta(t,x),$
 \be W_\delta(t,x)\leq \mbox{esssup}_{u_1
\in {\mathcal{U}}_{t,
t+\delta}}G^{t,x;u_1,\beta_1(u_1)}_{t,t+\delta} [W(t+\delta,
X^{t,x;u_1,\beta_1(u_1)}_{t+\delta})],\ \mbox{P-a.s.} \ee We use the
notation $I_\delta(t, x, u, v):=G^{t,x;u,v}_{t,t+\delta}
[W(t+\delta, X^{t,x;u,v}_{t+\delta})]$\ and notice that there exists
a sequence $\{u_i^1,\ i\geq 1\}\subset {\mathcal{U}}_{t, t+\delta}$\
such that
$$I_\delta(t, x, \beta_1):=\mbox{esssup}_{u_1 \in {\mathcal{U}}_{t,
t+\delta}}I_\delta(t, x, u_1, \beta_1(u_1))=\mbox{sup}_{i\geq
1}I_\delta(t, x, u_i^1, \beta_1(u_i^1)),\ \ \mbox{P-a.s.}$$ For any
$\varepsilon>0,$\ we put $\widetilde{\Gamma}_i:=\{I_\delta(t, x,
\beta_1)\leq I_\delta(t, x, u_i^1, \beta_1(u_i^1))+\varepsilon\}\in
{\mathcal{F}}_{t},\ i\geq 1.$\ Then
$\Gamma_1:=\widetilde{\Gamma}_1,\
\Gamma_i:=\widetilde{\Gamma}_i\backslash(\cup^{i-1}_{l=1}\widetilde{\Gamma}_l)\in
{\mathcal{F}}_{t},\ i\geq 2,$\ form an $(\Omega,
{\mathcal{F}}_{t})$-partition, and $u^\varepsilon_1:=\sum_{i\geq
1}\textbf{1}_{\Gamma_i}u_i^1$\ belongs obviously to
${\mathcal{U}}_{t, t+\delta}.$\ Moreover, from the nonanticipativity
of $\beta_1$\ we have $\beta_1(u^\varepsilon_1)=\sum_{i\geq
1}\textbf{1}_{\Gamma_i}\beta_1(u_i^1),$\ and from the uniqueness of
the solution of SDE (3.1) and RBSDE (3.5), we deduce that
$I_\delta(t, x, u^\varepsilon_1,
\beta_1(u^\varepsilon_1))=\sum_{i\geq
1}\textbf{1}_{\Gamma_i}I_\delta(t, x, u_i^1, \beta_1(u_i^1)),\
\mbox{P-a.s.}$\ Hence, \be
\begin{array}{llll}
W_\delta(t,x)\leq I_\delta(t, x, \beta_1)&\leq &\sum_{i\geq
1}\textbf{1}_{\Gamma_i}I_\delta(t, x, u_i^1, \beta_1(u_i^1))
+\varepsilon=I_\delta(t, x, u^\varepsilon_1,
\beta_1(u^\varepsilon_1))+\varepsilon\\
&=& G^{t,x;u^\varepsilon_1, \beta_1(u^\varepsilon_1)}_{t,t+\delta}
[W(t+\delta, X^{t,x;u^\varepsilon_1,
\beta_1(u^\varepsilon_1)}_{t+\delta})]+\varepsilon,\ \mbox{P-a.s.}
\end{array}
\ee
 On the other hand, using the fact that $\beta_1(\cdot):=\beta(\cdot\oplus u_2
)\in {\mathcal{B}}_{t, t+\delta}$\ does not depend on $u_2(\cdot)\in
{\mathcal{U}}_{t+\delta, T}$\ we can define
$\beta_2(u_2):=\beta(u^\varepsilon_1\oplus u_2)|_{[t+\delta, T]},\
\mbox{for all }\ u_2(\cdot)\in {\mathcal{U}}_{t+\delta, T}. $\ The
such defined $\beta_2: {\mathcal{U}}_{t+\delta, T}\rightarrow
{\mathcal{V}}_{t+\delta, T}$\ belongs to ${\mathcal{B}}_{t+\delta,
T}\ \mbox{since}\ \beta\in {\mathcal{B}}_{t, T}$. Therefore, from
the definition of $W(t+\delta,y)$\ we have, for any $y\in
{\mathbb{R}}^n,$
$$W(t+\delta,y)\leq \mbox{esssup}_{u_2 \in {\mathcal{U}}_{t+\delta, T}}J(t+\delta, y; u_2, \beta_2(u_2)),\ \mbox{P-a.s.}$$
Finally, because there exists a constant $C\in {\mathbb{R}}$\ such
that \be
\begin{array}{llll}
{\rm(i)} & |W(t+\delta,y)-W(t+\delta,y')| \leq C|y-y'|,\ \mbox{for any}\ y,\ y' \in {\mathbb{R}}^n;  \\
{\rm(ii)} & |J(t+\delta, y, u_2, \beta_2(u_2))-J(t+\delta, y',
u_2, \beta_2(u_2))| \leq C|y-y'|,\ \mbox{P-a.s.,}\\
 &\mbox{ }\hskip1cm \mbox{for any}\ u_2\in {\mathcal{U}}_{t+\delta, T},
\end{array}
\ee (see Lemma 3.2-(i) and (3.6)-(i)) we can show by approximating
$X^{t,x;u_1^\varepsilon,\beta_1(u_1^\varepsilon)}_{t+\delta}$\ by finite-valued ${\cal F}_{t+\delta}$-measurable random vectors that
$$W(t+\delta, X^{t,x;u_1^\varepsilon,\beta_1(u_1^\varepsilon)}_{t+\delta} )\leq
\mbox{esssup}_{u_2 \in {\mathcal{U}}_{t+\delta, T}}J(t+\delta,
X^{t,x;u_1^\varepsilon,\beta_1(u_1^\varepsilon)}_{t+\delta}; u_2,
\beta_2(u_2)),\ \mbox{P-a.s.}$$ To estimate the right side of the
latter inequality we note that there exists some sequence $\{u_j^2,\
j\geq 1\}\subset {\mathcal{U}}_{t+\delta, T}$\ such that
$$\mbox{esssup}_{u_2 \in {\mathcal{U}}_{t+\delta,
T}}J(t+\delta,X^{t,x;u_1^\varepsilon,\beta_1(u_1^\varepsilon)}_{t+\delta};
u_2, \beta_2(u_2))=\mbox{sup}_{j\geq
1}J(t+\delta,X^{t,x;u_1^\varepsilon,\beta_1(u_1^\varepsilon)}_{t+\delta};
u^2_j, \beta_2(u^2_j)),\ \mbox{P-a.s.}$$
 Then, putting\\
$\widetilde{\Delta}_j:=\{\mbox{esssup}_{u_2 \in
{\mathcal{U}}_{t+\delta,
T}}J(t+\delta,X^{t,x;u_1^\varepsilon,\beta_1(u_1^\varepsilon)}_{t+\delta};
u_2, \beta_2(u_2))\leq
J(t+\delta,X^{t,x;u_1^\varepsilon,\beta_1(u_1^\varepsilon)}_{t+\delta};
u^2_j, \beta_2(u^2_j))+\varepsilon\}\in {\mathcal{F}}_{t+\delta},\
j\geq 1;$\ we have with $\Delta_1:=\widetilde{\Delta}_1,\
\Delta_j:=\widetilde{\Delta}_j\backslash(\cup^{j-1}_{l=1}\widetilde{\Delta}_l)\in
{\mathcal{F}}_{t+\delta},\ j\geq 2,$\ an $(\Omega,
{\mathcal{F}}_{t+\delta})$-partition and
$u^\varepsilon_2:=\sum_{j\geq 1}\textbf{1}_{\Delta_j}u_j^2$\
 $\in {\mathcal{U}}_{t+\delta, T}.$ From
the nonanticipativity of $\beta_2$\ we have
$\beta_2(u^\varepsilon_2)=\sum_{j\geq
1}\textbf{1}_{\Delta_j}\beta_2(u_j^2),$\ and from the definition of
$\beta_1$ and $\beta_2$\ we know that $\beta(u_1^\varepsilon\oplus
u_2^\varepsilon)=\beta_1(u_1^\varepsilon)\oplus
\beta_2(u_2^\varepsilon ).$\ Thus, again from the uniqueness of the
solution of our FBSDE, we get
$$\begin{array}{lcl}
J(t+\delta,X^{t,x;u_1^\varepsilon,\beta_1(u_1^\varepsilon)}_{t+\delta};
u_2^\varepsilon,
\beta_2(u_2^\varepsilon))&=&Y_{t+\delta}^{t+\delta,X^{t,x;u_1^\varepsilon,\beta_1(u_1^\varepsilon)}_{t+\delta};
u_2^\varepsilon, \beta_2(u_2^\varepsilon)}\ \hskip2cm \mbox{(see (3.8))}\\
&=&\sum_{j\geq
1}\textbf{1}_{\Delta_j}Y_{t+\delta}^{t+\delta,X^{t,x;u_1^\varepsilon,\beta_1(u_1^\varepsilon)}_{t+\delta};
 u_j^2, \beta_2( u_j^2)}\\
&=&\sum_{j\geq
1}\textbf{1}_{\Delta_j}J(t+\delta,X^{t,x;u_1^\varepsilon,\beta_1(u_1^\varepsilon)}_{t+\delta};
u_j^2, \beta_2(u_j^2)),\ \mbox{P-a.s.}
\end{array}$$
Consequently, \be
\begin{array}{lll}
W(t+\delta,
X^{t,x;u_1^\varepsilon,\beta_1(u_1^\varepsilon)}_{t+\delta} )&\leq &
\mbox{esssup}_{u_2 \in {\mathcal{U}}_{t+\delta,
T}}J(t+\delta,X^{t,x;u_1^\varepsilon,\beta_1(u_1^\varepsilon)}_{t+\delta};
u_2, \beta_2(u_2))\\
&\leq & J(t+\delta,X^{t,x;u_1^\varepsilon,\beta_1(u_1^\varepsilon)}_{t+\delta};
u_2^\varepsilon,\beta_2(u_2^\varepsilon))+\varepsilon\\
&=& Y_{t+\delta}^{t+\delta,X^{t,x;u_1^\varepsilon,\beta_1(u_1^\varepsilon)}_{t+\delta};
u_2^\varepsilon,\beta_2(u_2^\varepsilon)}+\varepsilon\\
& = & Y_{t+\delta}^{t,x;u_1^\varepsilon\oplus u^\varepsilon_2,
\beta(u_1^\varepsilon\oplus
u^\varepsilon_2)}+\varepsilon\\
&=&Y_{t+\delta}^{t,x;u^\varepsilon,\beta(u^\varepsilon)}+\varepsilon,\
 \mbox{P-a.s.,}
\end{array}
\ee where $u^\varepsilon:= u_1^\varepsilon\oplus u^\varepsilon_2\in
{\mathcal{U}}_{t, T}.$\ From (6.21) and the definition of $W$ we conclude:
$$h(t+\delta, X^{t,x;u^\varepsilon,\beta(u^\varepsilon)}_{t+\delta})\le W(t+\delta,X^{t,x;u_1^\varepsilon,\beta_1(u_1^\varepsilon)}_{t+\delta})\leq
\left(Y^{t,x;u^\varepsilon,\beta(u^\varepsilon)}_{t+\delta}+\varepsilon\right)\wedge
h'(t+\delta, X^{t,x;u^\varepsilon,\beta(u^\varepsilon)}_{t+\delta}),
$$ $P$-a.s. From (6.19), (6.21) and the comparison theorem
for BSDEs with two reflecting boundaries we then get:
$$W_\delta(t,x)\leq G^{t,x;u_1^\varepsilon,\beta_1(u_1^\varepsilon)}_{t,t+\delta}
\left[\left(Y^{t,x;u^\varepsilon,\beta(u^\varepsilon)}_{t+\delta}+\varepsilon\right)\wedge
h'(t+\delta,
X^{t,x;u^\varepsilon,\beta(u^\varepsilon)}_{t+\delta})\right]+\varepsilon.$$
Thus, taking into account that
$Y^{t,x;u^\varepsilon,\beta(u^\varepsilon)}_{t+\delta}\le
h'(t+\delta, X^{t,x;u^\varepsilon,\beta(u^\varepsilon)}_{t+\delta}),$
we deduce from Lemma 2.9:

$W_\delta(t,x)\leq
G^{t,x;u_1^\varepsilon,\beta_1(u_1^\varepsilon)}_{t,t+\delta}
\left[\left(Y^{t,x;u^\varepsilon,\beta(u^\varepsilon)}_{t+\delta}+\varepsilon\right)\wedge
h'(t+\delta,
X^{t,x;u^\varepsilon,\beta(u^\varepsilon)}_{t+\delta})\right]+\varepsilon$

\qquad $\le
G^{t,x;u_1^\varepsilon,\beta_1(u_1^\varepsilon)}_{t,t+\delta}
\left[Y^{t,x;u^\varepsilon,\beta(u^\varepsilon)}_{t+\delta}\right]+(C+1)\varepsilon$

\qquad $= G^{t,x;u^\varepsilon,\beta(u^\varepsilon)}_{t,t+\delta}
\left[Y^{t,x;u^\varepsilon,\beta(u^\varepsilon)}_{t+\delta}\right]+(C+1)\varepsilon$

\qquad $=
Y^{t,x;u^\varepsilon,\beta(u^\varepsilon)}_{t}+(C+1)\varepsilon$

\qquad $\le $ esssup$_{u\in{\cal
U}_{t,T}}Y^{t,x;u,\beta(u)}_{t}+(C+1)\varepsilon$, $P$-a.s.

\smallskip
\noindent Therefore, in virtue of the arbitrariness of $\beta\in{\cal B}_{t,T},$
 \be W_\delta(t,x)\leq \mbox{essinf}_{\beta\in
{\mathcal{B}}_{t, T}}\mbox{esssup}_{u \in {\mathcal{U}}_{t,
T}}Y_{t}^{t,x;u,\beta(u)}+ (C+1)\varepsilon= W(t, x)+
(C+1)\varepsilon.\ee Finally, letting $\varepsilon\downarrow0,\
\mbox{we get}\ W_\delta(t,x)\leq W(t, x).$\endpf

\bl$ W(t, x)\leq W_\delta(t,x).$\el

 \noindent \textbf{Proof}. We continue to use the notations introduced above. From the definition of
$W_\delta(t,x)$\ we have
$$
\begin{array}{lll}
W_\delta(t,x)&=& \mbox{essinf}_{\beta_1 \in {\mathcal{B}}_{t,
t+\delta}}\mbox{esssup}_{u_1 \in {\mathcal{U}}_{t,
t+\delta}}G^{t,x;u_1,\beta_1(u_1)}_{t,t+\delta} [W(t+\delta,
X^{t,x;u_1,\beta_1(u_1)}_{t+\delta})]\\
&=&\mbox{essinf}_{\beta_1 \in {\mathcal{B}}_{t,
t+\delta}}I_\delta(t, x, \beta_1),
\end{array}
$$
where we have put
$$I_\delta(t, x, \beta_1):=\mbox{ esssup}_{u_1\in{\cal U}_{t,t+\delta}}I_\delta(t, x,u_1, \beta_1(u_1)).$$
We select $\{\beta_i^1,\ i\geq 1\}\subset{\mathcal{B}}_{t, t+\delta}$ such that
$$W_\delta(t,x)=\mbox{inf}_{i\geq 1}I_\delta(t, x, \beta_i^1),\ \mbox{P-a.s.,}$$ and for an arbitrarily small
$\varepsilon>0$\ we put $\widetilde{\Lambda}_i:=\{I_\delta(t, x,
\beta_i^1)-\varepsilon\leq W_\delta(t,x)\}\in {\mathcal{F}}_{t},\
i\geq 1,$ $\Lambda_1:=\widetilde{\Lambda}_1\ \mbox{and}\
\Lambda_i:=\widetilde{\Lambda}_i\backslash(\cup^{i-1}_{l=1}\widetilde{\Lambda}_l)\in
{\mathcal{F}}_{t},\ i\geq 2.$\ Then $\{\Lambda_i,\ i\geq 1\}$\ is an
$(\Omega, {\mathcal{F}}_{t})$-partition,
$\beta^\varepsilon_1:=\sum_{i\geq
1}\textbf{1}_{\Lambda_i}\beta_i^1$\ belongs to ${\mathcal{B}}_{t,
t+\delta},$\ and from the uniqueness of the solution of our FBSDE we
conclude that $I_\delta(t, x, u_1,
\beta^\varepsilon_1(u_1))=\sum_{i\geq
1}\textbf{1}_{\Lambda_i}I_\delta(t, x, u_1, \beta_i^1(u_1)),\
\mbox{P-a.s., for all}$\ \ $u_1(\cdot)\in {\mathcal{U}}_{t,
t+\delta}.$\ Hence,
 \be
\begin{array}{lll}
W_\delta(t,x)&\geq &\sum_{i\geq 1}\textbf{1}_{\Lambda_i}I_\delta(t,
x,
\beta_i^1)-\varepsilon\\
&\geq&\sum_{i\geq 1}\textbf{1}_{\Lambda_i}I_\delta(t, x, u_1,
\beta_i^1(u_1))-\varepsilon \\
&=& I_\delta(t, x, u_1, \beta^\varepsilon_1(u_1))-\varepsilon\\
&=& G^{t,x;u_1, \beta^\varepsilon_1(u_1)}_{t,t+\delta} [W(t+\delta,
X^{t,x;u_1, \beta_1^\varepsilon(u_1)}_{t+\delta})]-\varepsilon,\
\mbox{P-a.s., for all}\ \ u_1\in {\mathcal{U}}_{t, t+\delta}.
\end{array}
\ee
 On the other hand, from the definition of $W(t+\delta,y),$\
with the same technique as before, we deduce that, for any $y\in
{\mathbb{R}}^n,$\ there exists $\beta^\varepsilon_y\in
{\mathcal{B}}_{t+\delta, T}$\ \ such that \be W(t+\delta,y)\geq
\mbox{esssup}_{u_2 \in {\mathcal{U}}_{t+\delta,T}}J(t+\delta, y;
u_2, \beta^\varepsilon_y(u_2))-\varepsilon,\ \mbox{P-a.s.}\ee For every $m\ge 1$ we
now introduce
$$\varepsilon_m:=\left(\frac13\inf_{m-1\le\vert x\vert\le
m}(h'(t+\delta,x)-h(t+\delta,x))\right)\wedge \varepsilon >0.$$
Letting $C_0>0$ be the common Lipschitz constant of $h'(t+\delta,.)$
and $h(t+\delta,.)$ we put
$\delta_m=C_{0}^{-1}(\varepsilon_m\wedge\varepsilon)$. Moreover, we
let ${\cal O}_i^m,i\ge 1,$ be a Borel measurable decomposition of
$\Lambda_m:=\{x\in R^n:\, m-1\leq\vert x\vert<m\}$ such that
$\sum_{i\ge 1}{\cal O}_i^m=\Lambda_m$ and diam$({\cal O}_i^m)\leq
\delta_m,\, i\geq 1$. For each $m,i\ge 1$ we choose  an arbitrary
element of $y_i^m\in{\cal O}_i^m.$ Then, defining
$\left[X_{t+\delta}^{t,x;u_1,\beta_1^\varepsilon(u_1)}\right]:=
\sum_{i,m\ge
1}y_i^m\mathbf{1}_{\{X_{t+\delta}^{t,x;u_1,\beta_1^\varepsilon(u_1)}\in
{\cal O}_i^m\}}$, we have

\be\left\vert X_{t+\delta}^{t,x;u_1,
\beta_1^\varepsilon(u_1)}-\left[X_{t+\delta}^{t,x;u_1,
\beta_1^\varepsilon(u_1)}\right]\right\vert\leq\delta_m,\ee

\noindent and
$$\left\vert h\left(t+\delta,X_{t+\delta}^{t,x;u_1,\beta_1^\varepsilon(u_1)} \right)-h\left(t+\delta,\left[X_{t+\delta}^{t,x;u_1,\beta_1^\varepsilon(u_1)}\right]\right) \right\vert\leq\varepsilon_m\wedge\varepsilon,$$
everywhere on $\left\{X_{t+\delta}^{t,x;u_1,
\beta_1^\varepsilon(u_1)}\in\Lambda_m\right\}$, for all $u_1\in
{\cal U}_{t,t+\delta}.$ The same result also holds for $h'$ at the place
of $h$. Then,
$$h\left(t+\delta,X_{t+\delta}^{t,x;u_1,\beta_1^\varepsilon(u_1)}\right)\leq
h\left(t+\delta,\left[X_{t+\delta}^{t,x;u_1,\beta_1^\varepsilon(u_1)}\right]\right)+\varepsilon_m\hskip 4.2cm$$
$$\hskip 4.2cm  <h'\left(t+\delta,\left[X_{t+\delta}^{t,x;u_1,\beta_1^\varepsilon(u_1)}\right]\right)-\varepsilon_m\leq
h'\left(t+\delta,X_{t+\delta}^{t,x;u_1,\beta_1^\varepsilon(u_1)}\right),$$
on
$\left\{X_{t+\delta}^{t,x;u_1,\beta_1^\varepsilon(u_1)}\in\Lambda_m\right\},$
for all $u_1\in {\cal U}_{t,t+\delta}.$\ We choose for every $y_i^m$
some $\beta^{\varepsilon}_{y_i^m}\in{\cal B}_{t+\delta,T}$ such that
(6.24) is fulfilled, and clearly $\beta^{\varepsilon}_{u_1}:=
\sum_{i,m\ge 1}\mathbf{1}_{\{X_{t+\delta}^{t,x;u_1,
\beta_1^\varepsilon(u_1)}\in {\cal
O}_i^m\}}\beta^{\varepsilon}_{y_i^m}\in{\cal B}_{t+\delta,T}.$

 Now we can define the new strategy
$\beta^{\varepsilon}(u):=\beta_1^\varepsilon(u_1)\oplus
\beta^{\varepsilon}_{u_1}(u_2),\ u\in {\mathcal{U}}_{t, T},\
\mbox{where}\ u_1=u|_{[t, t+\delta]},\ u_2=u|_{(t+\delta, T]}$\
(restriction of $u$ to $[t, t+\delta]\times \Omega$\ and $(t+\delta,
T]\times \Omega$, resp.). Obviously, $\beta^{\varepsilon}$\ maps
${\mathcal{U}}_{t,T}$\ into ${\mathcal{V}}_{t,T}.$\ Moreover,
$\beta^{\varepsilon}$\ is nonanticipating: Indeed, let $S:
\Omega\longrightarrow[t, T]$\ be an $\mathbb{F}$-stopping time
and $u, u'\in {\mathcal{U}}_{t,T}$\ be such that $u\equiv u'$\ on
$\textbf{[\![}t, S\textbf{]\!]}$. Decomposing $u,\ u'$\ into $u_1,
u'_1\in {\mathcal{U}}_{t,t+\delta},\ u_2, u'_2\in
{\mathcal{U}}_{t+\delta, T}$\ such that $u=u_1\oplus u_2\
\mbox{and}\ u'=u'_1\oplus u'_2$\ we have $u_1\equiv u_1'$\ on
$\textbf{[\![}t, S\wedge(t+\delta)\textbf{]\!]}$,\ from where we get
$\beta_1^\varepsilon(u_1)\equiv \beta_1^\varepsilon(u_1')$\ on
$\textbf{[\![}t, S\wedge(t+\delta)\textbf{]\!]}$\ (recall that
$\beta_1^\varepsilon$\ is nonanticipating). On the other hand,
$u_2\equiv u_2'$\ on $\textbf{]\!]}t+\delta,
S\vee(t+\delta)\textbf{]\!]}(\subset (t+\delta,T]\times
\{S>t+\delta\}),$\ and on $\{S>t+\delta\}$\ we have $X^{t,x;u_1,
\beta_1^\varepsilon(u_1)}_{t+\delta}=X^{t,x;u'_1,
\beta_1^\varepsilon(u'_1)}_{t+\delta}.$\ Consequently, from our
definition, $\beta^{\varepsilon}_{u_1}=\beta^{\varepsilon}_{u'_1}$\
on $\{S>t+\delta\}$\ and
$\beta^{\varepsilon}_{u_1}(u_2)\equiv\beta^{\varepsilon}_{u'_1}(u'_2)$\
on $\textbf{]\!]}t+\delta, S\vee(t+\delta)\textbf{]\!]}.$ This
yields $\beta^{\varepsilon}(u)=\beta_1^\varepsilon(u_1)\oplus
\beta^{\varepsilon}_{u_1}(u_2)\equiv\beta_1^\varepsilon(u'_1)\oplus
\beta^{\varepsilon}_{u'_1}(u'_2)=\beta^{\varepsilon}(u')$\ on
$\textbf{[\![}t, S\textbf{]\!]}$, from where it follows that
$\beta^{\varepsilon}\in {\mathcal{B}}_{t, T}.$

Let now $u\in {\mathcal{U}}_{t, T}$\ be arbitrarily chosen and
decomposed into $u_1=u|_{[t, t+\delta]}\in {\mathcal{U}}_{t,
t+\delta}$\ and $u_2=u|_{(t+\delta, T]}\in {\mathcal{U}}_{t+\delta,
T}.$\ Then, with the notations

\smallskip

$W_{\varepsilon}(t+\delta,x):=\max\{h(t+\delta,x)+\varepsilon_m,\min\{W(t+\delta,x),h'(t+\delta,x)-\varepsilon_m\}\},$

$J_{\varepsilon}(t+\delta,x;u,v):=\max\{h(t+\delta,x)+\varepsilon_m,\min\{J(t+\delta,x;u,v),h'(t+\delta,x)-
\varepsilon_m\}\}$ and

$\widehat{J}_{\varepsilon}(t+\delta,x;u,v):=\max\{h(t+\delta,x)+\varepsilon_m,\min\{J(t+\delta,x;u,v)-\varepsilon,
h'(t+\delta,x)- \varepsilon_m\}\},$

\smallskip

\noindent for $x\in\Lambda_m,\, m\ge 1,\, u\in{\cal U}_{t+\delta,T}$ and $v\in{\cal V}_{t+\delta,T},$ we have  obviously,

\smallskip

\noindent$\widehat{J}_\varepsilon\left(t+\delta,\left[X_{t+\delta}^{t,x;u_1,
\beta_1^\varepsilon(u_1)}\right];u,v\right),\, J_\varepsilon\left(t+\delta,\left[X_{t+\delta}^{t,x;u_1,
\beta_1^\varepsilon(u_1)}\right];u,v\right), \, W_\varepsilon\left(t+\delta,\left[X_{t+\delta}^{t,x;u_1,
\beta_1^\varepsilon(u_1)}\right];u,v\right)$

$\in\left[h\left(t+\delta,X_{t+\delta}^{t,x;u_1,\beta_1^\varepsilon(u_1)}\right),h'\left(t+\delta,
 X_{t+\delta}^{t,x;u_1,\beta_1^\varepsilon(u_1)}\right)\right],$

\noindent and

$\left\vert J_\varepsilon\left(t+\delta,\left[X_{t+\delta}^{t,x;u_1,
\beta_1^\varepsilon(u_1)}\right];u,v\right)
-J\left(t+\delta,X_{t+\delta}^{t,x;u_1,
\beta_1^\varepsilon(u_1)};u,v\right)\right\vert\le C\varepsilon,$

\qquad $\left\vert W_\varepsilon\left(t+\delta,\left[X_{t+\delta}^{t,x;u_1,
\beta_1^\varepsilon(u_1)}\right]\right)
-W\left(t+\delta,X_{t+\delta}^{t,x;u_1,
\beta_1^\varepsilon(u_1)}\right)\right\vert\le C\varepsilon,$

\medskip

\noindent where $C$ is a constant which is independent of
$\varepsilon$ and the control processes. Thus, from (6.23),

\medskip

$W_\delta(t,x)\geq G_{t,t+\delta}^{t,x;u_1,
\beta_1^\varepsilon(u_1)}\left[W(t+\delta,X_{t+\delta}^{t,x;u_1,
\beta_1^\varepsilon(u_1)})\right]-\varepsilon$

$\geq G_{t,t+\delta}^{t,x;u_1,
\beta_1^\varepsilon(u_1)}\left[W_\varepsilon\left(t+\delta,\left[X_{t+\delta}^{t,x;u_1,
\beta_1^\varepsilon(u_1)}\right]\right)\right]-C\varepsilon$

$= G_{t,t+\delta}^{t,x;u_1,
\beta_1^\varepsilon(u_1)}\left[\sum_{i,m\ge
1}\mathbf{1}_{\{X_{t+\delta}^{t,x;u_1, \beta_1^\varepsilon(u_1)}\in
{\cal
O}_i^m\}}W_\varepsilon\left(t+\delta,y_i^m\right)\right]-C\varepsilon$

\smallskip

\noindent and, since
$$\left\vert \widehat{J}_\varepsilon\left(t+\delta,y;u,v\right)
-J_{\varepsilon}\left(t+\delta,y;u,v\right)\right\vert\le
\varepsilon,\, \mbox{ for all } y\in\mathbb{R}^n, u\in{\cal
U}_{t+\delta,T},v\in{\cal V}_{t+\delta,T},$$ \noindent we obtain

$W_\delta(t,x)\geq G_{t,t+\delta}^{t,x;u_1,
\beta_1^\varepsilon(u_1)}\left[\sum_{i,m\ge
1}\mathbf{1}_{\{X_{t+\delta}^{t,x;u_1, \beta_1^\varepsilon(u_1)}\in
{\cal
O}_i^m\}}\widehat{J}_\varepsilon\left(t+\delta,y_i^m;u_2,\beta^{\varepsilon}_{y_i^m}(u_2)\right)\right]-C\varepsilon$

$\geq G_{t,t+\delta}^{t,x;u_1,
\beta_1^\varepsilon(u_1)}\left[\sum_{i,m\ge
1}\mathbf{1}_{\{X_{t+\delta}^{t,x;u_1, \beta_1^\varepsilon(u_1)}\in
{\cal
O}_i^m\}}J_\varepsilon\left(t+\delta,y_i^m;u_2,\beta^{\varepsilon}_{y_i^m}(u_2)\right)\right]-C\varepsilon$

$= G_{t,t+\delta}^{t,x;u_1,
\beta_1^\varepsilon(u_1)}\left[J_\varepsilon\left(t+\delta,\left[X_{t+\delta}^{t,x;u_1,
\beta_1^\varepsilon(u_1)}\right];u_2,\beta^{\varepsilon}_{u_1}(u_2)\right)\right]-C\varepsilon$

$\geq G_{t,t+\delta}^{t,x;u_1,
\beta_1^\varepsilon(u_1)}\left[J\left(t+\delta,X_{t+\delta}^{t,x;u_1,
\beta_1^\varepsilon(u_1)};u_2,\beta^{\varepsilon}_{u_1}(u_2)\right)\right]-C\varepsilon$

$= G_{t,t+\delta}^{t,x;u,
\beta^\varepsilon(u)}\left[Y_{t+\delta}^{t,x;u,
\beta^\varepsilon(u)}\right]-C\varepsilon$

$=Y_{t}^{t,x;u, \beta^\varepsilon(u)}-C\varepsilon,$ $P$-a.s., for
any $u\in{\cal U}_{t,T}.$

\smallskip

\noindent (Notice that the constant $C$\ may be different from line
to line). This allows to conclude that $W_\delta(t,x)\geq
W(t,x)-C\varepsilon,$ and we get the wished relation by letting
$\varepsilon\rightarrow 0$.\endpf

\vskip0.5cm

\vskip 0.2cm

\noindent{\bf Acknowledgment} Juan Li thanks Mingyu Xu for some
helpful discussions on RBSDEs with two barriers.\

\end{document}